\documentclass[12pt]{article}
  \usepackage{amsfonts}
   \def\R{\mathbb{R}}
   \def\N{\mathbb{N}}
   
   \def\1{{\rm I\mskip -10.5mu 1}}
   
   \def\e{{\varepsilon}}
   \def\D{{\nabla}}
   \def\lo{\mathop{\longrightarrow}}
   
   \def\vi{{\varphi}}

   \def\cA{{\cal A}}
   \def\cB{{\cal B}}
   \def\cC{{\cal C}}
   \def\cD{{\cal D}}
   
   \def\cF{{\cal F}}
   \def\cG{{\cal G}}
   \def\cH{{\cal H}}

   \def\cN{{\cal N}}
   \def\cO{{\cal O}}

   \def\const{\mathop{\rm const}\nolimits}
   
   \def\loc{\mathop{\rm loc}\nolimits}

   \def\no{\noindent}
   \def\proof{\mbox {{\underline {\sf Proof}} \hspace{2mm}}}
   \def\qed{{\hfill {\em q.e.d.}\\\vspace{1mm}}}

   \newcommand{\beq}{\begin{equation}}
   \newcommand{\eeq}{\end{equation}}

\newtheorem{df}{Definition}[section]
\newtheorem{prop}[df]{Proposition}
\newtheorem{lemma}[df]{Lemma}
\newtheorem{teo}[df]{Theorem}
\newtheorem{rem}[df]{Remark}

\newtheorem{cor}[df]{Corollary}

 \newcommand{\sezione}[1]{\section{#1}\setcounter{equation}{0}}

\begin{document}


   \title{Multiple positive bound states for critical
     Schr\"odinger-Poisson systems}
  \author{Giovanna Cerami\thanks{ Dipartimento di Matematica,
        Politecnico di Bari,  Via Amendola 126/B -  70126
        Bari, Italia.
        E-mail: {\sf giovanna.cerami@poliba.it}}
    \and
        Riccardo Molle\thanks{Dipartimento di
        Matematica, Universit\`a di Roma ``Tor Vergata'', Via
                della Ricerca
        Scientifica n$^o$ 1 - 00133 Roma. E-mail: {\sf
                  molle@mat.uniroma2.it}}
        }
\date{}

 \maketitle



{\small {\sc \noindent \ \ Abstract.} -
Using variational methods we prove some results about existence and
multiplicity of positive bound states of to the following
Schr\"odinger-Poisson system: 
$$
\left\{ \begin{array}{l}
\vspace{2mm}
-\Delta u+V(x)u+K(x)\phi(x)u=u^5\\
-\Delta \phi =K(x)u^2\qquad x\in\R^3
\end{array}\right.\quad\quad (SP)
$$
We remark that $(SP)$ exhibits a ``double'' lack of compactness because
of the unboundedness of $\R^3$ and the critical growth of the
nonlinear term and that in our assumptions ground state solutions of
$(SP)$ do not exist. 

\vspace{3mm}



\vspace{1mm}


}


\sezione{Introduction}


This paper deals with the question of finding solutions to systems of the type

$$
\left\{ \begin{array}{l}
\vspace{2mm}
ih \frac{\partial \xi}{\partial t}= - \frac{h^2}{2m}\Delta \xi + (V(x)
+E) \xi +K(x)\phi(x)\xi - f(x,|\xi|)\\ 
-\Delta \phi =K(x)\xi^2\qquad x\in\R^3
\end{array}\right.\qquad\qquad 
$$
where $h$ is the Planck constant, $i$ is the imaginary unit, m is a
positive constant, E is a real number, $\xi : \R^3 \times [0,T]
\rightarrow  \textbf{C}.$ 

Such equations have strongly attracted the researchers  attention
because their deep physical meaning: we just mention they appear in
Semiconductor Theory  and Quantum Mechanics models (see f.i.the
celebrated papers \cite{BBL}, \cite{CPLL}, \cite{BF2}, \cite{PLL},
\cite{PLLS} and the book \cite{MRS}) 

Our interest is focused in the search of standing waves, that is
solutions  $ \xi(x,t) = e^ {- \frac{iEt}{h}}u $ where $u$ is a real
function. In this case one is  led to study the existence of functions
$u$ satisfying  
$$
\left\{ \begin{array}{l}
\vspace{2mm}
- \frac{h^2}{2m} \Delta u+V(x)u+K(x)\phi(x)u= f(x,|u|)\\
-\Delta \phi =K(x)u^2\qquad x\in\R^3.
\end{array}\right.\qquad\qquad (sp)
$$
Such a system is usually known as Schr\"odinger-Poisson system because
first equation, which is a nonlinear stationary  Schr\"odinger
equation, is  coupled with a Poisson equation.    
This model has been introduced in \cite{BF1} to describe electrostatic
situations in which the interaction between an electrostatic field and
solitary waves has to be considered; the nonlinear term $f$ simulates
the interaction between many particles, while, by the effect of the
Poisson equation,  the potential is determined by the charge of the
wave function itself.

More precisely we are interested in the existence of positive,
physically meaningful solutions of  $(sp)$ when $ f(x,u)= u^p$ and
$p=5 $ is the \textit{critical exponent} with respect to the Sobolev
embedding, moreover  we look for solutions when $h^2\over 2m$ is a
constant, which without loss of generality can be assumed equal to 1.
Therefore, the problem we address becomes 
$$
\left\{ \begin{array}{l}
\vspace{2mm}
-\Delta u+V(x)u+K(x)\phi(x)u=u^5\\
-\Delta \phi =K(x)u^2\qquad x\in\R^3.
\end{array}\right.\qquad\qquad (SP)
$$

It is well known that, in studying  Schr\"odinger-Poisson systems,
even in subcritical cases, one has to face many difficulties: some
come from the coupling and appear in the potential,  some originate
from  the lack of compactness due to the invariance of $\R^3$ under
translations. In addition, in the critical case the invariance by
dilations of $\R^3$ has to be considered and make the things even
harder to handle. 

During last fifteen years system $(sp)$ has been widely investigated,
mainly considering nonlinearities having subcritical growth:   first
results  have been obtained for equations with coefficients $V(x)$
and $ K(x) $ constants   or  radially symmetric;  then, an increasing
number of papers has been devoted to cases in which no symmetry
assumptions are requested. 

Describing the interesting and various contributions given to the
study of the subcritical case in a exhaustive  way, without forgetting
something, would be an hard task, so we prefer to refer readers
interested in a rich bibliography to the papers \cite{A}, \cite{BF3}
and to the book \cite{BF4}. 

Less contributions has received the analysis of $(sp)$ in the critical
case. Between those appeared in the latest years and treating
equations with non constant coefficients we refer the readers to
\cite{14CRH},  \cite{6HLW}, \cite{16HZ}, \cite{13LG}, \cite{8LG},
 \cite{3Z}, \cite{9Z} and references therein.  
However, we remark that the researchers
attention has been mainly devoted to the question of existence and
multiplicity of \textit{semi-classical} solutions (i.e. the case in
which, in $(sp),$  ${h^2\over 2m} \to 0,$ as in \cite{16HZ}, \cite{6HLW},
\cite{8LG}). Furthermore, we also stress  that all the researches we
are aware, concerning the search either of semi-classical either of
classical solution (as in \cite{14CRH}, \cite{13LG}, \cite{3Z}, \cite{9Z}), are
carried out under assumptions which  allow to work in frameworks that
ensure the existence of ground state solutions. 

Here we consider  situations that must  be faced by more refined tools. 
We ask the potentials $V(x)$ and $K(x)$ satisfy :

$$
\left\{ \begin{array}{lc}
\vspace{2mm}
\lim_{|x|\to +\infty} V(x)=V_\infty\ge 0 & (i)
\\ \vspace{2mm}
V(x)\ge V_\infty\quad\forall x\in\R^3 &(ii)
\\ \vspace{2mm}
 (V-V_\infty)\in L^{3/2}(\R^3)\hspace{5mm}& (iii)
\end{array}\right.\qquad\qquad\qquad\qquad  (H_V)
$$

and

$$
\left\{ \begin{array}{lc}
\vspace{2mm}
\lim_{|x|\to +\infty} K(x)= 0 & (i)
\\ \vspace{2mm}
K(x)\ge 0 \quad\forall x\in\R^3,\ K\not\equiv 0\hspace{5mm} &(ii)
\\ \vspace{2mm}
K\in L^{2}(\R^3).& (iii)
\end{array}\right. \qquad\qquad \qquad\qquad (H_K)
$$

Indeed, as shown in section 2,  Proposition \ref{P7.2}, under the
above  assumptions  the existence of positive solutions cannot be
obtained by minimization methods and ground state solutions do not
exist. 
Similar topological situations related to $(SP)$ have been considered
in the subcritical case in \cite{CV}, \cite{CM}. Here the critical
growth of the nonlinear term makes more difficult the question. 

The results we obtain are stated in the following  theorems, where $S$
denotes the best Sobolev constant. 

First theorem is concerned with potentials vanishing at infinity:

\begin{teo}
	\label{V=0}
	Let $V_\infty=0.$ Let  $(H_V)$, $(H_K),$ and
	
	\beq
	\label{stimaprecisabis}
	\left(1+{|V|_{L^{3/2}} \over S}+{|K|_{L^{2}}^2
		\over   S^{3/2}}\right)^3
	\,
	\left(1+{3\over 4\, S}\, |V|_{L^{3/2}}\right)<2
	\eeq  be satisfied.
	
	Then $(SP)$ has at least a positive solution.
\end{teo}

Second theorem provides existence and multiplicity of positive
solutions when $V_\infty \neq 0$, namely:

\begin{teo}
	\label{T2}
	Let $V_\infty>0.$  Let $(H_V)$ and $(H_K)$ be satisfied. 
	
	Then a real number $\bar V>0$ exists such that if $V_\infty\in
        (0,\bar V),$ then $(SP)$ has at least a positive solution. 
	
	Moreover, if in addition to  $(H_V)$ and $(H_K),$ 
	
	\beq
	\label{stimaprecisa}
	\left(1+{|V-V_\infty|_{L^{3/2}}\over S}+{|K|_{L^{2}}^2
		\over   S^{3/2}}\right)^3
	\,
	\left(1+{3\over 4\, S}\, |V-V_\infty|_{L^{3/2}}\right)<2
	\eeq  
	holds, $\bar V>0$ can be found so that,  when $V_\infty\in
        (0,\bar V),$ $(SP)$ has at least two distinct positive
        solutions. 
\end{teo}

It is worth observing that  Theorem \ref{V=0}   generalizes to
Schr\"odinger-Poisson systems a well known result proven in
\cite{BC90} for nonlinear Schr\"odinger equations. Indeed, in
\cite{BC90} the existence of a positive solution to 
$$
\left\{ \begin{array}{l}
\vspace{2mm}
-\Delta u+V(x)u = |u|^{\frac{N+2}{N-2}} \qquad x\in\R^N\\
u(x) \to 0 \qquad  \mbox{as}\ \  x\to \infty
\end{array}\right. \qquad \qquad\qquad \qquad (SE)
$$
has been shown assuming $V_\infty = 0,$  $V(x) \in L^{N/2},$ and, in
addition, a  restriction on  $|V|_{L^{3/2}}$ quite analogous to  the
bound $ |V|_{L^{3/2}}  <  (2^4 -1) S$ one can deduce from
(\ref{stimaprecisabis}) setting $K = 0$. 

On the other hand, as far as we know,  no resuls concerning solutions
of $(SE)$ when $V_\infty > 0$ are available, while a Pohozaev type
inequality shows that when $V$ is a positive constant $(SE)$ has no
solutions.  Therefore, it seems interesting to remark that from
Theorem \ref{T2} it follows,  as corollary, a non trivial  existence
and multiplicity theorem for $(SE)$.  We state it explicitely,
assuming $N=3$ because  Schr\"odinger-Poisson  systems are here
considered in $\R^3,$  nevertheless, reading the paper and the proof
of  Theorem \ref{T2}, it is not difficult to understand that
everything also holds for any dimension $N\geq 3.$  

\begin{teo}
	Let $ N= 3, V_\infty>0.$  Let $(H_V)$ be satisfied.
	
	Then a real number $\bar V>0$ exists such that if $V_\infty\in
        (0,\bar V),$ then $(SE)$ has at least a positive solution.  
	
	Moreover, if in addition to  $(H_V)$ 
	$$ |V-V_\infty|_{L^{3/2}}  <  (2^4 -1) S$$
	holds, $\bar V>0$ can be found so that,  when $V_\infty\in (0,\bar V),$ $(SE)$ has at least two distinct positive solutions. 
\end{teo}	

Of course some of our arguments, mainly those related to the lack of
compactness question, take advantage of some ideas introduced in
\cite{BC90}. Hovewer we strongly  point out that in the present paper
we face  different situations, we need new delicate estimates
concerning the nonlocal term, the variational framework in which we
work  must be different, and  we use here very refined and more subtle
tools to control translations and concentrations of Palais-Smale
sequences. 

The paper is organized as follows: in section 2 the variational
framework is introduced, some useful facts are stated, the compactness
question is studied and the nonexistence of ground state solutions is
proved,  section 3 contains some basic, deep estimates, and in section
4 the  proof of Theorems \ref{V=0} and \ref{T2} is performed.


\sezione{ Variational framework, Compactness\\ study, Nonexistence result}


Hereafter we use the following notation

{\small

\begin{itemize}

\item
$\|\cdot\|$ denotes the norm in $H$, that is
$$
\begin{array}{lc}
\vspace{2mm}
\|u\|=\left(\int_{\R^3}(|\D u|^2+V_\infty u^2)dx \right)^{1/2}& \mbox{ when
}H=H^1(\R^3)\\
\|u\|=\left(\int_{\R^3}|\D u|^2dx\right)^{1/2} & \mbox{ when
}H=\cD^{1,2}(\R^3).
\end{array}
$$

\item 
$|u|_q$, $1\le q\le +\infty$ denotes the norm in the Lebesgue space
$L^q(\R^3)$, while the norm of $u$ in $L^q(\Omega)$,
$\Omega\subset\R^3$, is denoted by $|u|_{q,\Omega}$.

\item
$B_\rho(y)$, $\forall y\in\R^3$, denotes the open ball of radius
$\rho$ centered at $y$, $(\cdot | \cdot )_{\R^3}$ denotes the scalar
product in $\R^3$, and for any measurable set $\cO\subset\R^3$,
$|\cO|$ denotes its Lebeasgue measure.

\item
$S$ is the best Sobolev constant, that is 
\beq
\label{1*}
S=\inf_{u\in H^1(\R^3)}{\|u\|^2\over |u|^2_6}
 =\inf_{u\in \cD^{1,2}(\R^3)}{\|u\|^2\over |u|^2_6}.
\eeq

\end{itemize}

}

\vspace{2mm}

\no \textit{Throughout the paper we set}
$$W(x):=V(x)-V_\infty$$
\textit{moreover we assume} $V$ \textit{and} $K$ \textit{satisfy} $(H_V)$ \textit{and}  $(H_K)$ \textit{ respectively}.

\vspace{3mm}

It is well known (see f.i. \cite{R}, \cite{CM}) that $(SP)$ can be
transformed in 
a nonlinear Schr\"odinger equation with a non local term.
Indeed, the Poisson equation can be solved by using the Lax-Milgram
theorem, thus for all $u\in\cD^{1,2}(\R^3)$ a unique
$\phi_u\in\cD^{1,2}(\R^3)$, satisfying 
\beq
\label{2.1}
-\Delta \phi=K(x)u^2,
\eeq is
obtained.
Then, inserting $\phi_u$ into the first equation of (SP) one gets
$$
-\Delta u+V(x)u+K(x)\phi_u(x)u=u^5.\qquad\qquad (SP')
$$
$(SP')$ is variational and its solutions are the critical points of the
functional 
$$
I(u)=\frac{1}{2}\left(\int_{\R^3}|\D u|^2dx+\int_{\R^3}V(x)u^2dx\right)
+\frac{1}{4}\int_{\R^3}K(x)\phi_u(x)
u^2dx-\frac{1}{6}\int_{\mathbb R^3}u^{6}dx
$$
which is defined in the space $H$ where $H$ is either $H^1(\R^3)$ or
$\cD^{1,2}(\R^3)$ according to whether $V_\infty>0$ or $V_\infty=0$.

Let $\Phi:H\to\cD^{1,2}(\R^3)$ be the operator defined by 
$$
\Phi(u)=\phi_u.
$$

Next two propositions collect some properties of $\Phi$.

\begin{prop}
\label{P2.1}
\begin{enumerate}
\item[1)] $\Phi$ is continuous;
\item[2)] $\Phi$ maps bounded sets into bounded sets;
\item[3)] $\Phi(tu)=t^2\Phi(u)$ for all $t\in \mathbb R$;
\item[4)] the following representation formula holds
\begin{equation}\label{conv}
\Phi (u)={1\over 4\pi}\int_{\R^3}{K(y)\over |x-y|}\, u^2(y)\,dy=
 {1\over 4\pi}{1\over |x|} \ast {K\, u^2}.
\end{equation}
\end{enumerate}
\end{prop}

The proof of {\em 1)} and {\em 2)} can be found f.i. in \cite{R} or \cite{CV}
while properties 
{\em 3)} and {\em 4)} are straight consequence of the fact that $\Phi(u)$ solves (\ref{2.1}).

\begin{prop}
\label{P2.2}
Let $(u_n)_n$, $u_n\in\cD^{1,2}(\R^3)$ be such that
\beq
\label{3.1old}
u_n\rightharpoonup 0\qquad\mbox{  in }\ \cD^{1,2}(\R^3).
\eeq
Then, up to subsequences
\begin{itemize}
\item[$a)$]  $\Phi(u_n)\to 0$ in $\cD^{1, 2} (\R^3)$
\item[$b)$]
  $\int_{\R^3}K(x)\phi_{u_n}u_n^2dx\lo 0$
\item[$c)$]
  $\int_{\R^3}K(x)\phi_{u_n}u_n\varphi\,
  dx\to 0$\quad $\forall\varphi\in \cD^{1, 2} (\R^3)$.
\end{itemize}
\end{prop}

The proof of Proposition \ref{P2.2} can be carried out exactly as that
of Proposition 2.2 of \cite{CM} once stated the following regularity
result, which is also useful in the study of the compactness
question.

\begin{lemma}
\label{L3.1}
Let $(u_n)_n$ be as in Proposition \ref{P2.2}.
Then, up to subsequences,
\beq
\label{3.1}
u_n\lo 0\quad\mbox{ in } L_{\loc}^p(\R^3)\quad\forall p\in[2,6).
\eeq
\end{lemma}

\no \proof By the Hardy inequality
$$
\int_{\R^3} u_n^2(1+|x|)^{-2}dx\le c\int_{\R^3}|\D u_n|^2dx\le C.
$$
Let $\cA\subset\subset\R^3$ be arbitrarily chosen and let $r>0$ be such that
$\cA\subset B_r(0)$, thus 
$$
\int_\cA u_n^2dx\le(1+r)^2\int_{B_r(0)}{u_n^2\over (1+|x|)^2}\, dx\le\hat
C
$$
which, together with (\ref{3.1old}) gives $(u_n)_n$ bounded in
$H^1(\cA)$.
Hence, up to a subsequence, $(u_n)_n$ converges strongly to 0 in $L^2(\cA)$ by the Sobolev
embedding and (\ref{3.1old}).
Then, the claim follows by interpolation and Sobolev embedding,
because for all $p\in [2,6)$
$$
|u_n|_{p,\cA}\le |u_n|_{2,\cA}^\alpha|u_n|_{6,\cA}^{1-\alpha}
$$
where ${\alpha\over 2}+{1-\alpha\over 6}={1\over p}$.

\qed

It is not difficult to verify that the functional $I$ is bounded
neither from below, nor from above.
So, it is suitable to consider $I$ restricted to the Nehari natural
constraint:
$$
\cN:=\{u\in H\setminus\{0\}\ :\ I'(u)[u]=0\}
$$
and remark that we can write $I_{|_\cN}$ as

\beq
\label{4.2}
\begin{array}{rcl}
\vspace{3mm}I_{|_{\cN}}(u)&
= &
{\displaystyle \frac{1}{4}\int_{\R^3}(|\D u|^2+V(x)u^2)dx+\frac{1}{12}\, \int_{\R^3}u^{6}dx}
\\ 
&=&
{\displaystyle \frac{1}{3}\, \int_{\R^3}(|\D u|^2+V(x)u^2)dx+
\frac{1}{12}\, \int_{\mathbb 
  R^3}K(x)\phi_u(x)u^2dx}
\end{array}
\eeq
from which one at once deduces that $I$ is bounded from below on
$\cN$. 
Furthermore, for all $u\in H\setminus\{0\}$, there exists a unique
$t_u\in\R^+\setminus\{0\}$ such that $t_uu\in\cN$.
Indeed, $t_u$ satisfies
\begin{eqnarray}
\nonumber
0=I'(tu)[tu] & = &
t^2\left(\int_{\R^3}|\D u|^2dx+\int_{\R^3}V(x)u^2dx\right)
+t^4\int_{\R^3}K(x)\phi_u(x)u^2dx\\
& & -t^{6}\int_{\R^3}u^{6}dx
\label{4.1}
\end{eqnarray}
which is easily seen to have a unique positive solution.

The function $t_u u\in\cN$ is called the {\bf projection of $u$ on
  $\cN$} and we also point out that
$$
I(t_uu)=\max_{t>0}I(tu).
$$
Actually, more precise information is available on $\cN$ and
$I_{|_\cN}$ and it can be summarized in the following lemma, whose
proof can be found in \cite{CV}

\begin{lemma}
\label{L5.1}
\begin{itemize}
\item[1)] $\cN$ is a $C^1$ regular manifold diffeomorphic to the sphere of $H$;
\item[2)] $I$ is bounded from below on $\cN$ by a positive
constant;
\item[3)] $u$ is a free critical point of $I$ if and only if $u$ is a critical point of $I$ constrained on $\cN$.
\end{itemize}
\end{lemma}
 
For what follows it is also useful to introduce the ``problem at
infinity'' related to (SP)
$$
\left\{ \begin{array}{l}
	-\Delta u=u^5 \qquad x\in\R^3\\
	u\in \cD^{1,2}(\R^3).
\end{array}\right.\qquad\qquad\qquad\qquad (P_\infty)
$$

Combining the results of \cite{GNN,G,T}  the following statement can
be obtained. 

\begin{prop}
\label{P5.1}
Any positive solution of $(P_\infty)$ must be of the form
\beq
\label{5.1}
\Psi_{\sigma,y}(x):=
{1\over\sigma^{1/2}}\,
\Psi
\left({x-y\over\sigma}\right)=
{[3\sigma^2]^{1/4}\over [\sigma^2+|x-y|^2]^{1/2}}
\eeq
where
$$
\Psi(x)=3^{1/4}\, {1\over (1+|x|^2)^{1/2}}={\Psi^*(x)\over |\Psi^*|_6}
$$
and
$$
\Psi^*(x)={1\over (1+|x|^2)^{1/2}}
$$
is the unique minimizer for $S$, up to translations and scaling.
\end{prop}

Throughout the paper we denote by
$$
I_\infty:H\to\R
$$
the functional whose critical points are solutions of $(P_\infty)$,
that is
$$
I_\infty(u)={1\over 2}\int_{\R^3}|\D u|^2dx-{1\over 6}\int_{\R^3}u^6dx
$$
and by 
$$
\cN_\infty=\{u\in H\setminus\{0\} :\ I_\infty'(u)[u]=0\}.
$$
A straight computation shows
$$
I_\infty(\Psi_{\sigma,y})=\min_{\cN_\infty}I_\infty(u)={1\over 3}\,
S^{3/2}.
$$

\begin{rem}
{\em
It is worth pointing out that the existence of infinitely many
changing sign solutions to $(P_\infty)$ has been proved by Ding
\cite{D}; however, for any such solution $u\in\cN_\infty$,
$u=u^+-u^-$, $u^+\neq 0\neq u^-$, the estimate $I_\infty(u)\ge {2\over
  3}\, S^{3/2}$ can easily be shown (see f.i. \cite{CSS}).

\vspace{2ex}

For all $u\in H\setminus\{0\}$, there exists unique
$\tau_u\in\R^+\setminus\{0\}$ such that $\tau_uu\in\cN_\infty$,
we call $\tau_uu$ projection of $u$ on $\cN_\infty$.
}
\end{rem}

\begin{lemma}
\label{L7.1}
Let $u\in H\setminus\{0\}$ and let $\tau_u u$ and $t_uu$ be the
projections of $u$ on $\cN_\infty$ and $\cN$ respectively.
Then
\beq
\label{7.1}
\tau_u\le t_u.
\eeq
\end{lemma}

\no \proof By definition we have 
$$
\tau_u^4={\int_{\R^3}|\D u|^2dx\over |u|_6^6}=
{
t_u^4|u|_6^6-t_u^2\int_{\R^3} K(x)\phi_uu^2dx-\int_{\R^3}V(x)u^2
\over
|u|_6^6
}
\le t_u^4.
$$
\qed

\begin{prop}
\label{P7.2}
Set
\beq
\label{7.2}
\inf\{I(u)\ :\ u\in\cN\}=:m.
\eeq
Then
\beq
\label{7.2bis}
m={1\over 3}\, S^{3/2}
\eeq
and the minimization problem (\ref{7.2}) has no solution.
\end{prop}

\no\proof
Let $u\in\cN$ be arbitrarily chosen and let  $\tau_uu$ be its projection on
$\cN_\infty$,
\begin{eqnarray*}
\vspace{2mm}
I(u) & \ge & I(\tau_uu)=I_\infty(\tau_uu)+{1\over
  2}\int_{\R^3}W(x)(\tau_uu)^2dx
\\ \vspace{2mm}
& &+{1\over 4}\int_{\R^3}K(x)\phi_{\tau_u u}(x)(\tau_uu)^2dx\ge
I_\infty(\tau_uu)\ge{1\over3}\, S^{3/2}
\end{eqnarray*}
from which
$$
m\ge {1\over 3}\, S^{3/2}
$$
follows.
To show that the equality holds, let us consider the sequence
$$
\tilde\Psi_n(x)=\chi(|x|)\Psi_{{1\over n},0}(x)=
\chi(|x|)\, {3^{1/4}\cdot\sqrt{1/n}\over \left((1/n)^2+|x|^2\right)^{1/2}}
$$
where $\chi\in\cC^\infty_0([0,+\infty))$ is a nonnegative real
function such that $\chi(s)=1$ if $s\in[0,1/2]$, and $\chi(s)=0$ if
$s\ge 1$.
Well known computations (see f.i. \cite{BN}) give
\beq
I_\infty(\tilde\Psi_n(x))={1\over
  2}\int_{\R^3}|\D\tilde\Psi_n(x)|^2dx-{1\over
  6}\int_{\R^3}\tilde\Psi^6_n(x)\,dx
= {1\over 3}S^{3/2}+O(1/n),
\label{8.1}
\eeq
\beq
\label{8.2}
|\tilde\Psi_n(x)|_2^2=O(1/n),
\eeq
\beq
\label{8.3}
\int_{\R^3}|\D\tilde\Psi_n(x)|^2dx-\int_{\R^3}\tilde\Psi^6_n(x)\,dx=O(1/n).
\eeq
Thus, $\tau_{\tilde\Psi_n}=1+O(1/n)$.
Now, using {\em 2)} of Proposition \ref{P2.1}, we get
\begin{eqnarray*}
\vspace{2mm}
&&\hspace{-5mm}\int_{\R^3}K(x)\phi_{\tilde\Psi_n}(x)\tilde\Psi_n^2(x)\,dx=
\int_{B_{1}(0)}K(x)\phi_{\tilde\Psi_n}(x)\tilde\Psi^2_n(x)dx
\\ \vspace{2mm}
&& =  \int_{B_{1}(0)\setminus
  B_{1/\sqrt{n}}(0)}K(x)\phi_{\tilde\Psi_n}(x)\tilde\Psi^2_n(x)\, dx
+
\int_{B_{1/\sqrt{n}}(0)}K(x)\phi_{\tilde\Psi_n}(x)\tilde\Psi^2_n(x)\, dx
\\ \vspace{2mm}
&&\le  |K|_2|\phi_{\tilde\Psi_n}|_6|\tilde\Psi_n|^2_{6,B_{1}(0)\setminus
  B_{1/\sqrt{n}}(0)}+|K|_{2,B_{1/\sqrt{n}}(0)}|\phi_{\tilde\Psi_n}|_6|\tilde\Psi_n|^2_6
\\
&&\le  c_1|\tilde\Psi_n|^2_{6,B_{1}(0)\setminus
  B_{1/\sqrt{n}}(0)} +c_2|K|_{2,B_{1/\sqrt{n}}(0)}
\end{eqnarray*}
with $c_1,c_2>0$ not depending on $n$.
Thus, considering that 
$$
|\tilde \Psi_n|_{6,B_{1}(0)\setminus
  B_{1/\sqrt{n}}(0)}=o(1),\qquad |K|_{2,B_{1/\sqrt{n}}(0)}=o(1)
$$
we deduce
\beq
\label{9.1}
\int_{\R^3}K(x)\phi_{\tilde\Psi_n}(x)\tilde\Psi^2_n(x)\, dx=o(1).
\eeq
Furthermore, for all $\rho>0$
\begin{eqnarray*}
\vspace{2mm}
\int_{\R^3}W(x)\tilde\Psi^2_n(x)\,
dx&=&\int_{B_\rho(0)}W(x)\tilde\Psi^2_n(x)\, dx+
\int_{\R^3\setminus B_\rho(0)}W(x)\tilde\Psi^2_n(x)\, dx\\
&\le &
|\tilde\Psi_n|_6^2|W|_{3/2,B_\rho(0)}+|W|_{3/2}|\tilde\Psi_n|^{2}_{6,
\R^3\setminus   B_\rho(0)}
\end{eqnarray*}
and, in view of
$$
\lim_{n\to \infty}|\tilde\Psi_n|_{6,\R^3\setminus B_\rho(0)}=0,
\quad 
\lim_{n\to\infty} |\tilde\Psi_n|_6=\const
$$
we obtain for all $\rho>0$
$$
\int_{\R^3}W(x)\tilde\Psi^2_n(x)\, dx\le \const |W|_{3/2,B_\rho(0)}.
$$
Thus, from $\lim\limits_{\rho\to0}|W|_{3/2,B_\rho(0)}=0$ we get
\beq
\label{9.2}
\int_{\R^3}W(x) \tilde\Psi^2_n(x)\, dx=o(1).
\eeq
Therefore, from (\ref{8.1}), (\ref{8.2}), (\ref{8.3}), (\ref{9.1}) and
(\ref{9.2})
\beq
\label{9.3}
t_{\tilde\Psi_n}=1+o(1)
\eeq
follows.

Finally, setting $\hat\Psi_n(x)=t_{\tilde\Psi_n}\tilde\Psi_n(x)$ we
conclude that (\ref{7.2bis}) holds because $\hat\Psi_n\in\cN$ and, by
(\ref{8.1}), (\ref{8.2}), (\ref{8.3}), (\ref{9.1}), (\ref{9.2}),
(\ref{9.3})
\begin{eqnarray*}
\vspace{2mm}
\lim_{n\to\infty}I(\hat\Psi_n)
& = &
{1\over
  2}t^2_{\tilde\Psi_n}\left[\int_{\R^3}(|\D\tilde\Psi_n(x)|^2+V(x)\tilde\Psi_n^2(x))dx\right]
\\ \vspace{2mm}
& & +{1\over
  4}t^4_{\tilde\Psi_n}\int_{\R^3}K(x)\phi_{\tilde\Psi_n}\tilde\Psi_n^2(x)\,dx
-{1\over  6}t^6_{\tilde\Psi_n}\int_\R^3\tilde\Psi_n^6(x)\, dx
\\
& = &
{1\over 3}\, S^{3/2}.
\end{eqnarray*}

To show that $m={1\over 3}\, S^{3/2}$ is not achieved we argue by
contradiction and we assume $\bar u\in\cN$ exists so that $I(\bar
u)={1\over 3}\, S^{3/2}$.
Then, by using (\ref{4.2}), $(H_k)(ii)$, $(H_V)(ii)$ we have 
\begin{eqnarray}
\vspace{2mm} \nonumber
{1\over 3}\, S^{3/2} & = & I(\bar u)={1\over 3}\int_{\R^3}(|\D\bar
u|^2+V(x)\bar u^2)dx+{1\over 12}\int_{\R^3}K(x)\phi_{\bar u}(x)\bar
u^2dx
\\ \vspace{2mm}\label{10.2}
&\ge& {1\over 3}\int_{\R^3}|\D\bar u|^2dx \ge{1\over
  3}\int_{\R^3}|\D\tau_{\bar u}\bar u|^2dx 
 \ge {1\over 3}\, S^{3/2}
\end{eqnarray}
from which
\beq
\label{10.1}
\int_{\R^3} K(x)\phi_{\bar u}\bar u^2dx=0,\quad \int_{\R^3}V(x)\bar
u^2dx=0
\eeq
follow, that, if $V_\infty\neq 0$ gives at once a contradiction.

When $V_\infty=0$, we observe that (\ref{10.1}), (\ref{10.2}), (\ref{7.1}) imply
$\tau_{\bar u}=1$, so
$$
\bar u(x)=\Psi_{\sigma,y}(x)>0\qquad\forall x\in\R^3
$$
for some $\sigma\in\R$, $y\in\R^3$. 
Then, by $(H_k)(ii)$
$$
\int_{\R^3}K(x)\phi_{\bar u}\bar u^2(x)\, dx>0
$$
has to be true, contradicting (\ref{10.1}).

\qed

Problem $(SP)$ is affected by a lack of compactness due to the
unboundedness of $\R^3$ and to the critical exponent.
Next proposition gives a picture of the compactness situation
describing the Palais-Smale sequences behaviour.

\begin{prop}
\label{P12.1}
Let $(u_n)_n$ be a PS-sequence of $I_{|_\cN}$, i.e. 
$$
u_n\in\cN,\qquad I(u_n)\lo c
$$
$$
\D I_{|_\cN}(u_n)\lo 0.
$$
Then there exist a number $k\in\N$, $k$ sequences of points
$(y_n^j)_n$, $y_n^j\in\R^3$, $1\le j\le k$, $k$ sequences of positive
numbers $(\sigma_n^j)_n$, $1\le j\le k$, $(k+1)$ sequences of
functions $(u_n^j)_n$, $u_n^j\in\cD^{1,2}(\R^3)$, $0\le j\le k$, such
that, replacing $(u_n)_n$, if necessary, by a subsequence still
denoted by $(u_n)_n$
\begin{eqnarray*}
\vspace{2mm}
i)&  & u_n(x)=u^0_n(x)+\sum_{j=1}^k{1\over (\sigma_n^j)^{1/2}}\,
u^j_n\left({x-y_n^j\over \sigma_n^j}\right)
\\
ii) & & u_n^j\lo u^j\ \mbox{ strongly in }\cD^{1,2}(\R^3),\ 0\le
j\le k
\end{eqnarray*}
where $u^0$ is a solution of $(SP')$, $u^j$ are solutions of
$(P_\infty)$, and
\begin{eqnarray*}
iii) & & \mbox{ if }y_n^j\to\bar y_j\quad \mbox{ then
}\quad\left\{\begin{array}{l}
\vspace{1mm} V_\infty\neq 0\ \Rightarrow\ \sigma_n^j\to 0\\
V_\infty=0 \ \Rightarrow\ 
\mbox{ either }\sigma_n^j\to 0\ 
\mbox{ or }\sigma_n^j\to\infty
\end{array}\right.
\\ 
 & & \mbox{ if }|y_n^j|\to\infty\quad\mbox{ then
}\quad\left\{\begin{array}{l}
\vspace{1mm} V_\infty\neq 0\ \Rightarrow\ \sigma_n^j\to 0
\\
V_\infty=0 \ \Rightarrow\ 
\mbox{ either }\sigma_n^j\to \bar\sigma^j\in\R^+\
\mbox{ or }\sigma_n^j\to\infty.
\end{array}\right.
\end{eqnarray*}
Moreover as $n\to\infty$
$$
\|u_n\|^2=\sum_{j=0}^k\|u^j\|^2+o(1)
$$
$$
I(u_n)\lo I(u^0)+\sum_{j=1}^k I_\infty(u^j).
$$
\end{prop}

The proof of Proposition \ref{P12.1} can be carried out quite analogously to
the proof of Theorem 2.5 in \cite{BC90} taking advantage of
Proposition \ref{P2.2}, Lemma \ref{L3.1} and {\em 3)} of Lemma \ref{L5.1}. 

\begin{cor}
\label{C13.1}
Assume $(u_n)_n$ satisfies the assumptions of Proposition \ref{P12.1}
with $c\in\left({1\over 3}\, S^{3/2},{2\over 3}\, S^{3/2}\right)$,
then $(u_n)_n$ is relatively compact.
\end{cor}

\no \proof
It suffices to apply Proposition \ref{P12.1} considering that any
nontrivial solution $u$ of $(SP')$, by Proposition \ref{P7.2},
verifies $I(u)>{1\over 3}\, S^{3/2}$, that any positive solution of
$(P_\infty)$ has energy ${1\over 3}\, S^{3/2}$, and for any changing
sign solution $v$ of $(P_\infty)$ the relation $I(v)\ge {2\over 3}\,
S^{3/2}$ holds.

\qed


\sezione{Basic estimates}


Let us introduce a barycenter type map 
$\beta:H\setminus\{0\}\to\R^3$:
$$
\beta(u)={1\over |u|_6^6}\int_{\R^3}{x\over 1+|x|}\, u^6(x)\, dx
$$
and a kind of inertial momentum $\gamma:H\setminus\{0\}\to\R$ to
estimate the concentration of a function $u$ around its barycenter:
$$
\gamma(u)={1\over |u|_6^6}\int_{\R^3}\left|{x\over 1+|x|}\,
  -\beta(u)\right|
u^6(x)\, dx.
$$
$\beta$ and $\gamma$ are continuous functions and 
\beq
\label{14.3}
\beta(tu)=\beta(u),\quad \gamma(tu)=\gamma(u)\qquad \forall t\in\R
\quad \forall u\in H \setminus\{0\}.
\eeq

\begin{prop} 
\label{P14.1}
The inequality
\beq
\label{14.1}
\inf\left\{I(u)\ :\ u\in\cN,\ \beta(u)=0,\ \gamma(u)={1/  2}\right\}>{1\over 3}\, S^{3/2}.
\eeq
holds true.
\end{prop}

\no \proof By (\ref{7.2bis}), clearly
\beq
\label{14.2}
\inf\left\{I(u)\ :\ u\in\cN,\ \beta(u)=0,\ \gamma(u)={1/2}\right\}\ge{1\over 3}\, S^{3/2}.
\eeq
To prove (\ref{14.1}) we argue by contradiction and suppose the
equality holds true in (\ref{14.2}).
Thus a sequence $(u_n)_n$ exists such that
\beq
\label{15.1}
\left\{
\begin{array}{lcc}
\vspace{2mm}
u_n\in\cN,\ \beta(u_n)=0,\ \gamma(u_n)={1/2} & & (a) \\
\lim\limits_{n\to\infty}I(u_n)={1\over 3}\, S^{3/2}. & & (b)
\end{array}\right.
\eeq
Therefore, since $K(x)\ge 0$, $V(x)\ge 0$, $\phi_{u_n}\ge 0$, by using
(\ref{4.2}) and (\ref{7.1}), we can write
\begin{eqnarray}
\vspace{2mm} \nonumber
{1\over 3}\, S^{3/2} & = & \lim\limits_{n\to\infty}\left[{1\over
    3}\int_{\R^3}(|\D u_n|^2+V(x)u_n^2)dx+{1\over
    12}\int_{\R^3}K(x)\phi_{u_n}u_n^2\right] 
\\ \label{15.3} \vspace{2mm}
& \ge & \lim\limits_{n\to\infty}{1\over
    3}\int_{\R^3}|\D u_n|^2dx\ge  \lim\limits_{n\to\infty}{1\over
    3}\, \tau^2_{u_n}\int_{\R^3}|\D u_n|^2dx
\\  \nonumber
& \ge &{1\over 3}\, S^{3/2} ,
\end{eqnarray}
from which we infer 
$$
\lim_{n\to\infty}{1\over
    3}\int_{\R^3}|\D u_n|^2dx={1\over 3}\, S^{3/2}
$$
and
\beq
\label{15.2}
\lim_{n\to\infty}\tau_{u_n}=1.
\eeq
By the uniqueness of the family of ground state positive solutions of
$(P_\infty)$ stated in Proposition \ref{P5.1}, and by Proposition \ref{P12.1}
we deduce 
$$
\tau_{u_n}u_n(x)=\Psi_{\sigma_n,y_n}(x)+\e_n(x)
$$
where $\sigma_n\in\R$, $\sigma_n>0$, $y_n\in\R^3$,
$\e_n\in\cD^{1,2}(\R^3)$, $\e_n\to 0$ strongly in $\cD^{1,2}(\R^3)$
and $L^6(\R^3)$.
Furthermore, by (\ref{15.2})
\beq
\label{16.2}
u_n(x)=\Psi_{\sigma_n,y_n}(x)+\tilde\e_n(x)
\eeq
with $\tilde\e_n\to 0$ strongly in $\cD^{1,2}(\R^3)$ and $L^6(\R^3)$.

We claim that, up to subsequences,
\beq
\label{16.1}
a)\quad \lim_{n\to\infty}\sigma_n=\bar\sigma>0\qquad
b)\lim_{n\to\infty}y_n=\bar y\in\R^3.
\eeq
Indeed, once the claim is shown true, the proof can be quickly
concluded: it is enough to observe that in this case
\beq
\label{16.3}
\Psi_{\sigma_n, y_n}\to\Psi_{\bar\sigma,\bar y}\ \mbox{ strongly in
}\cD^{1,2}(\R^3)\ \mbox{ and }L^6(\R^3)
\eeq
so (\ref{15.1})$(b)$, (\ref{16.2}), and (\ref{16.3}), together with
$(H_K)(ii)$, $\Psi_{\bar\sigma,\bar y}>0$,
$\phi_{\Psi_{\bar\sigma,\bar y}}>0$ give
\begin{eqnarray*}
\vspace{2mm}
S^{3/2}&=& \lim\limits_{n\to\infty}\left[\int_{\R^3}(|\D
  u_n|^2+V(x)u_n^2)dx+{1\over
    4}\int_{\R^3}K(x)\phi_{u_n}u_n^2dx\right]
\\ \vspace{2mm} 
& = & \int_{\R^3}(|\D
  \Psi_{\bar\sigma,\bar y}|^2+V(x)\Psi_{\bar\sigma,\bar y})dx+{1\over
    4}\int_{\R^3}K(x)\phi_{\Psi_{\bar\sigma,\bar y}}\Psi_{\bar\sigma,\bar y}^2dx
\\
& > & \int_{\R^3}|\D
  \Psi_{\bar\sigma,\bar y}|^2=S^{3/2}
\end{eqnarray*}
that is impossible.

Let us now prove the claim. To show (\ref{16.1})$(a)$, we start
estabilishing  that $(\sigma_n)_n$  is bounded.
Assume $(\sigma_n)_n$ unbounded, then, passing
eventually to a subsequence, $\sigma_n\lo \infty$ occurs, thus for all
$\rho>0$
\beq
\label{17.1}
\lim_{n\to\infty}\int_{B_\rho(0)}u^6_n
dx=\lim_{n\to\infty}\int_{B_\rho(0)}\Psi_{\sigma_n,y_n}^6(x)\, dx=0.
\eeq
Hence, from (\ref{15.1})$(a)$ we get
\beq
\label{17.2}
0=\beta(u_n)=\beta(\Psi_{\sigma_n,y_n})+o(1),
\eeq
and, taking into account {\em 2)} of Lemma \ref{L5.1} we obtain for all $\rho>0$
\begin{eqnarray*}
\vspace{2mm}
\gamma(u_n) & = &{1\over |u_n|_6^6}\int_{\R^3}{|x|\over 1+|x|}\,
  u_n^6(x)\, dx
\\ \vspace{2mm}
& = &
{1\over |u_n|_6^6}\left[ \int_{B_\rho(0)}{|x|\over 1+|x|}\,
  u^6_n(x)\,dx+\int_{\R^3\setminus B_\rho(0)}{|x|\over 1+|x|}\,
  u^6_n(x)\, dx\right]
\\
& = &
{1\over |u_n|_{6,\R^3\setminus B_\rho(0)}^6+o(1)} 
\left[\int_{\R^3\setminus B_\rho(0)}{|x|\over 1+|x|}\,
  u^6_n(x)\, dx +o(1)\right]
\\ 
& \ge &
{\rho\over 1+\rho}+o(1),
\end{eqnarray*}
so $\liminf\limits_{n\to\infty}\gamma(u_n)\ge {\rho\over 1+\rho}$,
$\forall \rho >0$, which implies
$$
\lim_{n\to\infty}\gamma(u_n)=1
$$
contradicting (\ref{15.1})$(a)$.

Thus, up to a subsequence, $\sigma_n\lo\bar\sigma\in\R^+$.
If $\bar\sigma=0$ would occur, then for all $\rho>0$
$$
\lim_{n\to\infty}\int_{\R^3\setminus B_\rho(y_n)}u_n^6(x)\,dx=
\lim_{n\to\infty}\int_{\R^3\setminus B_\rho(y_n)}(\Psi_{\sigma_n,y_n}(x))^6dx=0
$$
and
$$
0<c<|u_n|_6^6=|u_n|_{6,B_\rho(y_n)}^6+o(1),
$$
from which, for all $\rho>0$
\begin{eqnarray*}
\vspace{2mm}
{|y_n|\over 1+|y_n|} & = & \left|{y_n\over
    1+|y_n|}-\beta(u_n)\right|={1\over
  |u_n|_6^6}\left|\int_{\R^3}\left({y_n\over 1+|y_n|}-{x\over 1+|x|}\right)u^6_n(x)\,dx\right|
\\ \vspace{2mm}
& \le & {1\over
  |u_n|_6^6}\left[\int_{B_\rho(y_n)}\left|
{y_n\over 1+|y_n|}-{x\over
      1+|x|}\right|u_n^6(x)\, dx \right.
\\ \vspace{2mm} 
& & {\phantom{{1\over
  |u_n|_6^6}**}}\left.
 + \int_{\R^3\setminus B_\rho(y_n)}\left|{y_n\over
      1+|y_n|}-{x\over 1+|x|}\right|u^6_n(x)\,dx\right]
\\ \vspace{2mm}
& \le &
{1\over |u_n|^6_{6,B_\rho(y_n)}+o(1)}\left[
  \rho|u_n|^6_{6,B_\rho(y_n)}+o(1)\right]
\\
&\le & \rho+o(1)
\end{eqnarray*}
which implies $|y_n|\to 0$ as $n\to \infty$.
Therefore 
\begin{eqnarray*}
0\le\gamma(u_n) & = & {1\over |u_n|_6^6}\int_{\R^3}\left|{x\over
    1+|x|}-{y_n\over 1+|y_n|}\right|\, u_n^6(x)\, dx+o(1)
\\
& \le &\rho+o(1)\qquad\forall \rho>0.
\end{eqnarray*}
So, we obtain $\lim\limits_{n\to\infty}\gamma(u_n)=0$ against
(\ref{15.1})$(a),$ therefore (\ref{16.1})$(a)$ is proved.

Let us now show that $(|y_n|)_n$ is bounded and, then, convergent up
to subsequences.
By contradiction we suppose that a subsequence, still denoted by
$(y_n)_n$, exists for which $\lim\limits_{n\to\infty}|y_n|=\infty$.
Then for all $\e>0$, and all $R>0$, $\bar n\in\N$ can be found so that
$\forall n>\bar n$
$$
|x-y_n|<R\ \Rightarrow\ \left|{x\over 1+|x|}-{y_n\over
    1+|y_n|}\right|<\e.
$$
Moreover, for all $\e>0$, $\bar\rho>0$ depending only
on $\e$  exists  such that $\forall R >\bar\rho$
\beq
\label{19.1}
\int_{\R^3\setminus B_R(y_n)}\Psi_{\bar\sigma,y_n}^6(x)\, dx=
\int_{\R^3\setminus B_R(0)}\Psi_{\bar\sigma,0}^6(x)\, dx <\e.
\eeq
Now, let us choose arbitrarily $\e>0$ and fix $R>0$ so that
(\ref{19.1}) holds true; for large $n$ we get
\begin{eqnarray*}
\vspace{2mm}
\left|\beta(u_n)-{y_n\over 1+|y_n|}\right| & \le &
{1\over |u_n|_6^6}\int_{\R^3}\left|{x\over 1+|x|}-{y_n\over
    1+|y_n|}\right|\, u_n^6(x)\, dx
\\ \vspace{2mm} & \le &
{1\over |\Psi_{\bar\sigma,y_n}|^6_6+o(1)}\left[
\int_{B_R(y_n)}\left|{x\over 1+|x|}-{y_n\over
    1+|y_n|}\right|\, \Psi_{\bar\sigma,y_n}^6(x)\, dx\right.
\\ \vspace{2mm}
& &\left. {\phantom{{1\over \Psi^6_{\bar\sigma,y_n}}}}
+\int_{\R^3\setminus B_R(y_n)}\left|{x\over 1+|x|}-{y_n\over
    1+|y_n|}\right|\, \Psi_{\bar\sigma,y_n}^6(x)\, dx\right]+o(1)
\\
& \le & \hat c\e
\end{eqnarray*}
with $\hat c>0$ independent of $y_n$ and $R$. Thus $|\beta(u_n)|\lo 1$
as $n\to\infty$, giving a contradiction with (\ref{15.1})$(a)$ and
completing the proof of the claim and of the proposition.

\qed

\begin{prop}
\label{P21.1}
Let assume $V_\infty > 0$. Set
$$
\mu=\inf\{I(u)\ :\ u\in\cN,\ \beta(u)=0,\ \gamma(u)\ge 1/2\}.
$$
Then
\beq
\label{21.3}
\mu>{1\over 3}\, S^{3/2}.
\eeq
\end{prop}

\no \proof We follow an analogous argument to that of Proposition
\ref{P14.1}.
We start observing that by (\ref{7.2bis})
\beq
\label{21.1}
\mu\ge{1\over 3}\, S^{3/2}
\eeq
and if the equality in (\ref{21.1}) holds a sequence $(u_n)_n$ exists
so that
\beq
\label{21.2}
\left\{\begin{array}{lcc}
\vspace{2mm}
u_n\in\cN,\ \beta(u_n)=0,\ \gamma(u_n)\ge 1/2 & & (a)\\
\lim\limits_{n\to\infty} I(u_n)={1\over 3}\, S^{3/2}. & & (b)
\end{array}\right.
\eeq
Then, the same computations made in Proposition \ref{P14.1} allow to
assert that
$$
u_n(x)=\Psi_{\sigma_n,y_n}(x)+\e_n(x)
$$
where $\sigma_n\in\R$, $\sigma_n>0$, $y_n\in\R^3$, $\e_n\to 0$ in
$\cD^{1,2}(\R^3)$.

The sequence $(\sigma_n)_n$ is bounded.
Indeed, otherwise, up to
subsequences, $\sigma_n\to \infty$ and, by (\ref{21.2}), $(H_V)(ii)$,
$(H_K)(ii)$,  (\ref{4.2}),
\begin{eqnarray*}
\vspace{2mm}
{1\over 3}\, S^{3/2} & = &
\lim\limits_{n\to\infty}I(u_n)\ge \liminf\limits_{n\to\infty}\left({1\over
  3}\int_{\R^3}|\D u_n|^2dx+V_\infty \int_{B_{\sigma_n}(y_n)}u_n^2dx\right)
\\ \vspace{2mm}
& \ge &
 \lim\limits_{n\to\infty}{1\over
  3}\int_{\R^3}|\D \Psi_{\sigma_n,y_n}|^2dx+
\lim\limits_{n\to\infty}V_\infty\sigma_n^2\left[\int_{B_1(0)} 
\Psi_{1,0}^2(x)\, dx+o(1)\right]
\\
& = & + \infty.
\end{eqnarray*}
Hence, passing eventually to a subsequence, the relation
$\lim\limits_{n\to\infty}\sigma_n=\bar\sigma$ holds.
Working again as in Proposition (\ref{14.1}), $\bar\sigma>0$ is shown
and, furthermore, $(y_n)_n$ bounded is proved, so that $y_n\to\bar y$,
up to a subsequence.
Thus, we deduce
$$
\Psi_{\sigma_n,y_n}\lo\Psi_{\bar\sigma,\bar y}\quad\mbox{ strongly in
}\cD^{1,2}(\R^3)\mbox{ and }L^2_{\loc}(\R^3)
$$
and the impossible relation
\begin{eqnarray*}
\vspace{2mm}
{1\over 3}\, S^{3/2} & = & \lim\limits_{n\to\infty} I(u_n)\ge 
{1\over 3}\left[\int_{\R^3} |\D\Psi_{\bar\sigma ,\bar
    y}|^2dx+V_\infty\int_{B_{\bar\sigma}(\bar y)}\Psi^2_{\bar\sigma,\bar
    y}dx\right]
\\
& > & 
{1\over 3}\int_{\R^3} |\D\Psi_{\bar\sigma,\bar
    y}|^2dx={1\over 3}\, S^{3/2}
\end{eqnarray*}
which brings to conclude the equality in (\ref{21.1}) cannot occur.

\qed

\begin{rem}
{\em
Let us remark that  if $V_\infty=0$ then $\mu={1\over 3}\, S^{3/2}$.

}\end{rem}

In what follows we use the notation
$$
\cB_{V_\infty}:=\inf\{I(u)\ :\ u\in\cN,\ \beta(u)=0,\
\gamma(u)={1/2}\},\qquad\mbox{ if }V_\infty>0,
$$
and
$$
\cB_{0}:=\inf\{I(u)\ :\ u\in\cN,\ \beta(u)=0,\
\gamma(u)={1/2}\},\qquad\mbox{ if }V_\infty=0.
$$
Obviously, by Proposition \ref{P14.1},
$$
\cB_0>{1\over 3}\, S^{3/2}\quad\mbox{ and }\quad \cB_{V_\infty}>{1\over 3}\, S^{3/2}.
$$
\begin{rem}
{\em
Let us observe that
\beq
\label{23.1}
\cB_0\le\cB_{V_\infty}\qquad\forall V_\infty\in\R^+\setminus\{0\}.
\eeq
}\end{rem}

Indeed, for all $u\in H^{1}(\R^3)$ such that $\beta(u)=0$,
$\gamma(u)={1/3}$ and
$$
\int_{\R^3}(|\D
u|^2+W(x)u^2)dx+\int_{\R^3}K(x)\phi_uu^2dx-\int_{\R^3}u^6dx=0
$$
let us consider $t_uu$ such that 
$$
t_u^2\int_{\R^3}[|\D u|^2+V(x)u^2]dx
+t_u^4\int_{\R^3}K(x)\phi_uu^2dx-t_u^6\int_{\R^3}u^6dx=0,
$$
then $\beta(t_uu)=\beta(u)=0$,
$\gamma(t_uu)=\gamma(u)={1/3}$ and
$$
\hspace{-1.5cm} 
{1\over 2}
t_u^2\int_{\R^3}(|\D u|^2+V(x)u^2)dx
+{1\over 4}t_u^4\int_{\R^3}K(x)\phi_uu^2dx-{1\over
  6}t_u^6\int_{\R^3}u^6dx
$$

\vspace{-6mm}

\begin{eqnarray*}
& \ge &
{1\over 2}
\int_{\R^3}(|\D u|^2+(V_\infty+W(x))u^2)\, dx
+{1\over 4}\int_{\R^3}K(x)\phi_uu^2dx-{1\over
  6}\int_{\R^3}u^6dx
\\ 
& \ge &
{1\over 2}
\int_{\R^3}(|\D u|^2+W(x)u^2)dx
+{1\over 4}\int_{\R^3}K(x)\phi_uu^2dx-{1\over
  6}\int_{\R^3}u^6dx
\end{eqnarray*}
from which considering also $\cD^{1,2}(\R^3)\supset H^{1}(\R^3)$
(\ref{23.1}) follows.

\vspace{3mm}

Let us now fix a number $\bar c$ such that 
\beq
\label{23'.0}
{1\over 3}\, S^{3/2}<\bar c<\min\left({\cB_0+{1\over 3}\, S^{3/2}\over
    2}, {1\over 2}\,  S^{3/2}\right).
\eeq
Moreover, if either (\ref{stimaprecisabis}) or (\ref{stimaprecisa})
holds, then we also require
\beq
\label{stimaprecisa2}
\bar c <{2\over 3}\ S^{3/2} \left(1+{|W|_{3/2}\over S}+
{|K|^2_2\over S^{3/2}}\right)^{-3}\,
\left(1+{3\over 4\, S}\, |W|_{3/2}\right)^{-1}
\eeq
 and we observe that $W(x) = V(x)$ when $V_\infty =0$.

 Clearly $\bar c\in \left({1\over 3}\, S^{3/2}, {2\over 3}\,
  S^{3/2}\right).$

We denote by $\omega(x)$ a function belonging to $H^1(\R^3)$ having
the following properties:
\beq
\label{23'*}
\left\{
\begin{array}{cl}
(i) &   \omega\in\cC_0^\infty(B_1(0)),\\
(ii) & \omega(x)\ge 0\qquad\forall x\in B_1(0)\\
(iii) & \omega\in\cN_\infty\mbox{ and
}I_\infty(\omega)=\Sigma\in\left({1\over 3}\,S^{3/2},\bar c\right)\\
(iv) & \omega(x)=\omega(|x|)\mbox{ and } |x_1|\le |x_2|\
\Rightarrow\ \omega(x_1)\ge\omega(x_2). 
\end{array}\right.
\eeq

Moreover, for every $\sigma>0$ and $y\in\R^3$ we set
$$
\omega_{\sigma,y}(x)\ =\ \left\{\begin{array}{lc}\vspace{2mm}
\sigma^{-1/2}\omega\left({x-y\over\sigma}\right) & x\in B_\sigma(y)\\
0 & x\not\in B_\sigma(y).
\end{array}\right.
$$
Remark that
\beq
\label{23'.1}
|\omega|_6=|\omega|_{6,B_1(0)}=|\omega_{\sigma,y}|_{6,B_\sigma(y)}=|\omega_{\sigma,y}|_6.
\eeq

\begin{lemma}
\label{L24.1}
 The following relations hold
\beq
\label{24.1}
\begin{array}{rcl}
\vspace{2mm}
a) & & \lim\limits_{\sigma\to 0}\sup\left\{\int_{\R^3}
  W(x)\, \omega^2_{\sigma,y}(x)\, dx\ :\ y\in\R^3\right\}=0\\
\vspace{2mm}
b) & &  \lim\limits_{\sigma\to \infty}\sup\left\{\int_{\R^3}
  W(x)\, \omega^2_{\sigma,y}(x)\, dx\ :\ y\in\R^3\right\}=0\\
\vspace{2mm}
c) & & \lim\limits_{r\to \infty}\sup\left\{\int_{\R^3}
  W(x)\, \omega^2_{\sigma,y}(x)\, dx\ :\ y\in\R^3,\ |y|=r,\ \sigma>0\right\}=0.
\end{array}
\eeq
\end{lemma}

\no \proof
Let $y\in\R^3$ be arbitrarily chosen. Then for all $\sigma>0$
$$
\int_{\R^3} W(x)\, \omega_{\sigma,y}^2(x)\, dx=
\int_{B_\sigma(y)} W(x)\, \omega_{\sigma,y}^2(x)\, dx\le
|W|_{3/2,B_\sigma(y)}|\omega|^2_{6,B_1(0)}\le c|W|_{3/2,B_\sigma(y)}
$$
with $c$ not depending on $\sigma$. So
$$
\sup_{y\in\R^3}\int_{\R^3} W(x)\omega_{\sigma,y}^2(x)\, dx\le c
\sup\{|W|_{3/2,B_\sigma(y)}\ ;\ y\in\R^3\}
$$
and, since $\lim\limits_{\sigma\to
  0}|W|_{3/2,B_\sigma(y)}=0$, uniformly in $y\in\R^3$,  (\ref{24.1})$(a)$ follows.

To prove  (\ref{24.1})$(b)$, again let us choose arbitrarily
$y\in\R^3$.
Then $\forall\rho>0$, $\forall \sigma>0$
\begin{eqnarray*}
\vspace{2mm}
\int_{\R^3}W(x)\, \omega^2_{\sigma,y}dx & = & 
\int_{B_\rho(0)}W(x)\, \omega^2_{\sigma,y}dx +
\int_{\R^3\setminus B_\rho(0)}W(x)\, \omega^2_{\sigma,y}dx
\\ \vspace{2mm}
&\le &
|W|_{3/2,B_\rho(0)}|\omega_{\sigma,y}|^2_{6,B_\rho(0)}+
|W|_{3/2,\R^3\setminus B_\rho(0)}|\omega_{\sigma,y}|^2_{6,\R^3\setminus
  B_\rho(0)}
\\
&\le &
|W|_{3/2,B_\rho(0)}\sup\limits_{y\in\R^3}|\omega_{\sigma,y}|^2_{6,B_\rho(0)}+
c|W|_{3/2,\R^3\setminus B_\rho(0)}
\end{eqnarray*}
with $c$ not depending on $\sigma$ and $\rho$.

Now, 
$\lim\limits_{\sigma\to \infty}|\omega_{\sigma,y}|_{6,B_\rho(0)}=0$,
uniformly with respect to $y\in\R^3$,
so we get $\forall\rho>0$
$$
\lim_{\sigma\to\infty}\sup_{y\in\R^3}\int_{\R^3}W(x)\, \omega^2_{\sigma,y}(x)\,
dx\le c|W|_{3/2,\R^3\setminus B_\rho(0)}
$$
and, letting $\rho\to +\infty$, (\ref{24.1})$(b)$ follows.

To verify (\ref{24.1})$(c)$, we argue by contradiction, so we assume
the existence of a sequence $(y_n)_n$, $y_n\in\R^3$, and a sequence
$(\sigma_n)_n$, $\sigma_n\in\R^+\setminus\{0\}$, such that 
\beq
\label{25.1}
|y_n|\lo \infty
\eeq
and
\beq
\label{25.2}
\lim_{n\to\infty}
\int_{\R^3}W(x)\, \omega^2_{\sigma_n,y_n}(x)\, dx >0.
\eeq
In view of (\ref{24.1})$(a)-(b)$, passing eventually to a subsequence,
we can suppose $\lim\limits_{n\to\infty}\sigma_n=\bar\sigma>0$ which,
together (\ref{25.1}) and $(H_v)\, (iii)$, implies
$$
\lim_{n\to\infty} |W|_{3/2,B_{\sigma_n(y_n)}}=0.
$$
Therefore we deduce the relation
$$
\lim_{n\to\infty}\int_{\R^3}W(x)\, \omega^2_{\sigma_n,y_n}(x)\, dx\le
\lim_{n\to\infty}\left[|W|_{3/2,B_{\sigma_n}(y_n)}\cdot
  |\omega_{\sigma_n,y_n}|^2_{6,B_{\sigma_n}(y_n)}\right]=0
$$
contradicting (\ref{25.2}).

\qed

\begin{lemma}
\label{L26.1}
The following relations hold
\beq
\label{26.1}
\begin{array}{rcl}
\vspace{2mm}
a) & & \lim\limits_{\sigma\to 0}\sup\left\{\int_{\R^3}
  K(x)\, \phi_{\omega_{\sigma,y}}(x)\, \omega^2_{\sigma,y}(x)\, dx\ :\ y\in\R^3\right\}=0\\
\vspace{2mm}
b) & &  \lim\limits_{\sigma\to \infty}\sup\left\{\int_{\R^3}
  K(x)\, \phi_{\omega_{\sigma,y}}(x)\, \omega^2_{\sigma,y}(x)\, dx\ :\ y\in\R^3\right\}=0\\
\vspace{2mm}
c) & & \lim\limits_{r\to \infty}\sup\left\{\int_{\R^3}
  K(x)\, \phi_{\omega_{\sigma,y}}(x)\, \omega^2_{\sigma,y}(x)\, dx\ :\ y\in\R^3,\ |y|=r,\ \sigma>0\right\}=0.
\end{array}
\eeq
\end{lemma}

\no\proof
We start remarking that, by {\em 2)} of Proposition \ref{P2.1},
$\{\phi_{\omega_{\sigma,y}}\ :\ \sigma>0,\ y\in\R^3\}$ is a bounded
set in $\cD^{1,2}(\R^3)$ and $L^6(\R^3)$.
Let $y\in\R^3$ arbitrarily chosen, then for all $\sigma>0$
\begin{eqnarray*}
\vspace{2mm}
\int_{\R^3}K(x)\phi_{\omega_{\sigma,y}}(x)\,  \omega^2_{\sigma,y}(x)\,
dx & = & \int_{B_\sigma(y)}K(x)\phi_{\omega_{\sigma,y}}(x)\,
\omega^2_{\sigma,y}(x)\, 
dx 
\\ \vspace{2mm}
& \le &
|K|_{2,B_\sigma(y)}|\phi_{\omega_{\sigma,y}}|_{6}|
\omega_{\sigma,y}|^2_{6,B_\sigma(y)} 
\\
& \le &
c|K|_{2,B_\sigma(y)}
\end{eqnarray*}
$c>0$ independent of $y$ and $\sigma$.
Hence
$$
\sup_{y\in\R^3}\int_{\R^3}K(x)\phi_{\omega_{\sigma,y}}(x)\,  \omega^2_{\sigma,y}(x)\,
dx\le c\sup_{y\in\R^3}|K|_{2,B_\sigma(y)}
$$
which gives (\ref{26.1})$(a)$ because
$$
\lim_{\sigma\to 0} |K|_{2,B_\sigma(y)}=0\qquad\mbox{ uniformly
  in } y\in\R^3.
$$
To verify (\ref{26.1})$(b)$, let us fix arbitrarily $y\in\R^3$.
Then $\forall \rho>0$ $\forall\sigma>0$
\begin{eqnarray*}
\vspace{2mm}
 & &\hspace{-2cm}\int_{\R^3}K(x)\phi_{\omega_{\sigma,y}}(x)\,  \omega^2_{\sigma,y}(x)\,
dx
\\ \vspace{2mm}
& = & \
\int_{\R^3\setminus B_\rho(0)}K(x)\phi_{\omega_{\sigma,y}}(x)\,  \omega^2_{\sigma,y}(x)\,
dx
+
\int_{ B_\rho(0)}K(x)\phi_{\omega_{\sigma,y}}(x)\,  \omega^2_{\sigma,y}(x)\,
dx
\\ \vspace{2cm}
& \le & |K|_{2,\R^3\setminus B_\rho(0)}|\phi_{\omega_{\sigma,y}}|_6|\omega_{\sigma,y} |^2_{6,\R^3\setminus B_\rho(0)}
+|K|_{2, B_\rho(0)}|\phi_{\omega_{\sigma,y}}|_6|\omega_{\sigma,y} |^2_{6,B_\rho(0)}
\\ 
& \le &
\bar c_1|K|_{2,\R^3\setminus B_\rho(0)}|+\bar c_2 |K|_{2,
  B_\rho(0)}|\sup\limits_{y\in\R^3}| \omega_{\sigma,y} |^2_{6,B_\rho(0)}
\end{eqnarray*}
with $\bar c_1,\bar c_2\in\R^+\setminus\{0\}$ depending neither on  $y$
nor on $\sigma$.
Thus, considering $|\omega_{\sigma,y}|_{6,B_\rho(0)}\to 0$ as
$\sigma\to\infty$, uniformly with respect to $y\in\R^3$, we get
$$
\lim_{\sigma\to\infty}\sup\left\{\int_{\R^3}K(x)\phi_{\omega_{\sigma,y}}(x)\,  \omega^2_{\sigma,y}(x)\,
dx \ :\ y\in\R^3\right\}\le\bar c_1|K|_{2,\R^3\setminus B_\rho(0)}
$$
and, letting $\rho\to\infty$, we obtain (\ref{26.1})$(b)$.
To show (\ref{26.1}), working by contradiction, we assume it false, so
that $(y_n)_n$, $y_n\in\R^3$ and $(\sigma_n)_n$,
$\sigma_n\in\R^3\setminus\{0\}$, exist for which 
\beq
\label{27.1}
(a)\quad |y_n|\to\infty, \qquad (b)\quad \lim_{n\to\infty} 
\int_{\R^3}K(x)\phi_{\omega_{\sigma_n,y_n}}(x)\,  \omega^2_{\sigma_n,y_n}(x)\,
dx >0.
\eeq
By (\ref{26.1})$(a)-(b)$, up to a subsequence, we have
$\lim\limits_{n\to\infty}\sigma_n=\bar\sigma\in(0,\infty)$, and using
(\ref{27.1})$(a)$ we deduce
$$
\lim_{n\to\infty}| K|_{2,B_{\sigma_n}(y_n)}=0.
$$
Thus
\begin{eqnarray*}
\int_{\R^3}K(x)\phi_{\omega_{\sigma_n,y_n}}(x)
\omega^2_{\sigma_n,y_n}(x)\, dx &= &
\int_{B_{\sigma_n}(y_n)}K(x)\phi_{\omega_{\sigma_n,y_n}}(x)\omega^2_{\sigma_n,y_n}(x)\,
dx \\
& \le &
|K|_{2,B_{\sigma_n}(y_n)}|\phi_{\omega_{\sigma_n,y_n}}|_6|\omega_{\sigma_n,y_n}|^2\le\tilde
c |K|_{2,B_{\sigma_n}(y_n)}
\end{eqnarray*}
$\tilde c$ not depending on $y_n$ nor on $\sigma_n$. Therefore,
$$
\lim_{n\to\infty} 
\int_{\R^3}K(x)\phi_{\omega_{\sigma_n,y_n}}(x)
\omega^2_{\sigma_n,y_n}(x)\, dx \le \tilde c
\lim_{n\to\infty} 
|K|_{2,B_{\sigma_n}(y_n)}=0
$$
follows and contradicts (\ref{27.1})$(b)$.

\qed

\begin{cor}
\label{C28.1}
Set $t_{\sigma,y}=t_{\omega_{\sigma,y}}$.
Then
\beq
\label{28.1}
\begin{array}{rcl}
\vspace{2mm}
a) & & \lim\limits_{\sigma\to 0}\sup\left\{|t_{\sigma,y}-1|\ :\ y\in\R^3\right\}=0\\
\vspace{2mm}
b) & &  \lim\limits_{\sigma\to \infty}\sup\left\{|t_{\sigma,y}-1|\  :\ y\in\R^3\right\}=0\\
\vspace{2mm}
c) & & \lim\limits_{r\to \infty}\sup\left\{|t_{\sigma,y}-1|\
:\ y\in\R^3,\ |y|=r,\ \sigma>0\right\}=0.
\end{array}
\eeq
\end{cor}

\no \proof It is enough to observe that, since
$\omega_{\sigma,y}\in\cN_\infty$, the equalities 
$$
1={\|\omega_{\sigma,y}\|^2\over |\omega_{\sigma,y}|^6_6}= 
{
t_{\sigma,y}^4|\omega_{\sigma,y}|^6_6-t_{\sigma,y}^2
\int_{\R^3}K\phi_{\omega_{\sigma,y}}  \omega^2_{\sigma,y}-
\int_{\R^3}W\omega^2_{\sigma,y}
\over|\omega_{\sigma,y}|^6_6
}
$$
hold true, so by Lemmas \ref{L24.1} and \ref{L26.1}, relations
(\ref{28.1}) follow at once.

\qed

\begin{lemma}
\label{L29.1}
The following relations hold
\beq
\label{29.1}
\begin{array}{rcl}
\vspace{2mm}
a) & & \lim\limits_{\sigma\to 0}\sup\left\{\gamma(\omega_{\sigma,y})\ :\ y\in\R^3\right\}=0\\
\vspace{2mm}
b) & &  \lim\limits_{\sigma\to
  \infty}\inf\left\{\gamma(\omega_{\sigma,y})\  :\ y\in\R^3,\ |y|\le
  r\right\}=1  \quad \forall r>0\\
\vspace{2mm}
c) & & (\beta(\omega_{\sigma,y})|y)_{\R^3}>0\quad\forall y\in\R^3,\ \forall\sigma>0.
\end{array}
\eeq
\end{lemma}

\no \proof For all $y\in\R^3$ and for all $\sigma>0$
\begin{eqnarray*}
0\le\gamma(\omega_{\sigma,y}) & = &
{1\over |
  \omega_{\sigma,y}|^6_{6,B_\sigma(y)}}\int_{B_\sigma(y)}\left|{x\over
      1+|x|}-\beta(\omega_{\sigma,y})\right|\omega^6_{\sigma,y}\, dx
  \\
& \le & 
{1\over |
  \omega_{\sigma,y}|^6_{6,B_\sigma(y)}}\int_{B_\sigma(y)}\left|{x\over
      1+|x|}-{y\over
      1+|y|}\right|\omega^6_{\sigma,y}\, dx+
\left|{y\over
      1+|y|}- \beta(\omega_{\sigma,y})\right|\\
& \le & 
{1\over |
  \omega_{\sigma,y}|^6_{6,B_\sigma(y)}}\left[\int_{B_\sigma(y)}|x-y|\,
  \omega^6_{\sigma,y}dx+\int_{B_\sigma(y)}\left|{x\over
      1+|x|}-{y\over
      1+|y|}\right|\omega^6_{\sigma,y}\, dx\right]\\
& \le & 2\sigma
\end{eqnarray*}
from which
$$
0\le\sup\{\gamma(\omega_{\sigma,y})\ :\ y\in\R^3\}\le 2\sigma
$$
follows and, then, (\ref{29.1})$(a)$.

In order to prove (\ref{29.1})$(b)$, let us first show that for all
$y\in\R^3$
\beq
\label{30.1}
\lim_{\sigma\to\infty}\beta(\omega_{\sigma,y})=0.
\eeq
Indeed, by symmetry $\beta(\omega_{\sigma,0})=0$, $\forall\sigma>0$ 
considering (\ref{23'.1}), we deduce
\begin{eqnarray*}
|\beta(\omega_{\sigma,y})| & = &
{1\over |
  \omega_{\sigma,y}|^6_{6}}\left|\int_{\R^3}{x\over
      1+|x|}\, \omega^6_{\sigma,y} dx\right|
\\
& = &
{1\over |
  \omega_{\sigma,0}|^6_{6}}\left|\int_{\R^3}{x\over
      1+|x|}(\omega^6_{\sigma,y} -\omega^6_{\sigma,0}) dx\right|
\\
& \le &
{1\over |
  \omega_{\sigma,0}|^6_{6}}\, \int_{\R^3}{|x|\over
      1+|x|}|\omega^6_{\sigma,y} -\omega^6_{\sigma,0}| dx
\\
& \le & c\int_{\R^3}|\omega^6_{1,y/\sigma} -\omega^6_{1,0}| dx
\lo 0\quad\mbox{ as }\sigma\to\infty.
\end{eqnarray*}
Now, let us choose arbitrarily $r>0$ and $y\in\R^3$ such that $|y|\le
r$.
For all $\sigma>0$, we have
$$
\gamma(\omega_{\sigma,y})={1\over |\omega_{\sigma,y}|^6_6}\int_{\R^3}
\left|{x\over
      1+|x|}-\beta(\omega_{\sigma,y})\right|\omega^6_{\sigma,y}\,  
dx\le 1+|\beta(\omega_{\sigma,y})|
$$
from which, by (\ref{30.1})
$$
\limsup_{\sigma\to\infty} \inf\{\gamma(\omega_{\sigma,y})\ :\
y\in\R^3,\ |y|\le r\}\le 1
$$
follows.
Hence, to obtain (\ref{29.1})$(b)$ we just need to show 
\beq
\label{new}
\liminf_{\sigma\to\infty} \inf\{\gamma(\omega_{\sigma,y})\ :\
y\in\R^3,\ |y|\le r\}\ge 1.
\eeq
To do this let us take the sequences $(y_n)_n$,   $y_n\in\R^3$, $
|y_n|\le r$,  and $(\sigma_n)_n,$  $ \sigma_n \in \R^+,$  
 $\sigma_n\to\infty$.
Then, considering (\ref{30.1}) and that $|\omega_{\sigma,y}|_{6,B_\rho(0)}\to0$, as
$\sigma\to\infty$, we deduce for all $\rho>0$
\begin{eqnarray*}
\gamma(\omega_{\sigma_n,y_n}) & = &
{1\over |\omega_{\sigma_n,y_n}|^6_6}\int_{\R^3}
\left|{x\over
      1+|x|}-\beta(\omega_{\sigma_n,y_n})\right|\omega^6_{\sigma_n,y_n} dx
\\
& \ge & 
{1\over |\omega_{\sigma_n,y_n}|^6_{6,\R^3\setminus B_\rho(0)}+o(1)}\int_{\R^3\setminus B_\rho(0)}
{|x|\over
      1+|x|}\, \omega^6_{\sigma_n,y_n}dx-o(1) 
\\
& \ge & 
{\rho\over 1+\rho}-o(1)
\end{eqnarray*}
which, letting $\rho\to\infty$, gives
$\lim\limits_{n\to\infty}\gamma(\omega_{\sigma_n,y_n})=1$ and proves (\ref{new}).

Lastly, let us remark that (\ref{29.1})$(c)$ is immediate if $0\not\in
B_\sigma(y)$.
If $0\in B_\sigma(y)$, to prove  (\ref{29.1})$(c)$ we just need to
consider that $\forall \bar x\in B_\sigma(y)$ such that $(\bar x|y)_{\R^3}>0$
the point  $-\bar x$ verifies  $\omega_{\sigma,y}(-\bar x)<
\omega_{\sigma,y}(\bar x)$.

\qed


\sezione{Proof of Theorems}


In this section we use the notation $I_0$ to denote the functional $I$
when $V_\infty=0$, while when we write $I$ we intend that $V_\infty\neq 0$.
Namely, $\forall u\in\cD^{1,2}(\R^3)$
$$
I_0(u)={1\over 2}\int_{\R^3} (|\D u|^2+W(x)u^2)dx+{1\over
    4}\int_{\R^3}K(x)\phi_u u^2dx-{1\over 6}\int_{\R^3}u^6dx
$$
and, $\forall u\in H^{1}(\R^3)$
$$
I(u)={1\over 2}\int_{\R^3} (|\D u|^2+(V_\infty+W(x))u^2)dx+{1\over
    4}\int_{\R^3}K(x)\phi_u u^2dx-{1\over 6}\int_{\R^3}u^6dx
$$
with $V_\infty\neq 0$.

According to this we have
$$
\cN_0:=\{u\in\cD^{1,2}(\R^3)\setminus\{0\}\ :\ I'_0(u)[u]=0\}
$$
$$
\cN:=\{u\in H^{1}(\R^3)\setminus\{0\}\ :\ I'(u)[u]=0\}
$$
and we denote by
$$
\hat\omega_{\sigma,y}:=t_{\sigma,y,0}\omega_{\sigma,y}=t^0_{\omega_{\sigma,y}}\omega_{\sigma,y}
$$
$$
\tilde\omega_{\sigma,y}:=t_{\sigma,y}\omega_{\sigma,y}=t_{\omega_{\sigma,y}}\omega_{\sigma,y}
$$
the projections of $\omega_{\sigma,y}$ respectively on $\cN_0$ and on
$\cN$.

\begin{lemma}
\label{L34.1}
There exist real numbers $\bar r>0$, $\sigma_1,\sigma_2:\
0<\sigma_1<{1\over2}<\sigma_2$ such that
\beq
\label{34.1}
\gamma(\hat\omega_{\sigma_1,y})<1/2,\
\gamma(\hat\omega_{\sigma_2,y})>1/2,
\quad\forall y\in\R^3
\eeq
and
\beq
\label{34.2}
\sup\{I_0(\hat\omega_{\sigma,y})\ :\ (\sigma,y)\in\partial\cH\}<\bar c
\eeq
where $\bar c$ is defined in (\ref{23'.0}) and
\beq
\cH=\{(\sigma,y)\in\R^+\times\R^3\ :\ \sigma\in[\sigma_1,\sigma_2],\
|y|<\bar r\}.
\eeq
\end{lemma}

\no\proof
Since
\beq
\label{34.4}
\begin{array}{rcl}
\vspace{2mm}
I_0(\hat\omega_{\sigma,y})&=&
{1\over 2}t^2_{\sigma,y,0}\int_{\R^3}|\D\omega_{\sigma,y}|^2
-{1\over  6}t^6_{\sigma,y,0}\int_{\R^3}\omega_{\sigma,y}^6dx
\\ & & +{1\over 2}t^2_{\sigma,y,0}\int_{\R^3}W(x)\omega_{\sigma,y}^2dx
+{1\over 4}t^4_{\sigma,y,0}\int_{\R^3}K(x)\phi_{\omega_{\sigma,y}}
\omega_{\sigma,y}^2dx
\end{array}
\eeq
the existence of $\sigma_1\in\left(0,{1/2}\right)$ such that
$\gamma(\hat\omega_{\sigma,y})<{1/2}$ and $I_0(\hat\omega_{\sigma,y}) < \bar{c}$ holds true when
$\sigma=\sigma_1$ for all $y\in\R^3$ is a consequence of (\ref{14.3}),
(\ref{24.1})$(a)$, (\ref{26.1})$(a)$, (\ref{28.1})$(a)$,
(\ref{29.1})$(a)$, and (\ref{23'*}).
Furthermore (\ref{24.1})$(c)$, (\ref{26.1})$(c)$, (\ref{28.1})$(c)$ and
(\ref{23'*}) allow us to choose $\bar r>0$ such that, if $|y|=\bar r,\ $
$I_0(\hat\omega_{\sigma,y}) < \bar{c}$ is satisfied for all $\sigma>0$. 
Once $\bar r$ is fixed, we use again (\ref{14.3}), (\ref{23'*}) plus 
(\ref{24.1})$(b)$, (\ref{26.1})$(b)$, (\ref{28.1})$(b)$,
(\ref{29.1})$(b)$ in (\ref{34.4}) and we find $\sigma_2>{1/2}$
for which
$$
\gamma(\hat\omega_{\sigma_2,y})>1/2\quad\forall y\in\R^3,\ |y|\le\bar
r
$$
and $I_0(\hat\omega_{\sigma,y}) < \bar{c}$  is verified when $\sigma=\sigma_2$ and $|y|\le\bar
r$.

\qed

\begin{lemma}
\label{L36}
Let $\sigma_1,\sigma_2,\bar r,\cH$ as in Lemma \ref{L34.1}.
Then, $(\tilde\sigma,\tilde y)\in\partial\cH$ and $(\bar\sigma,\bar
y)\in\stackrel{\circ}{\cH}$ exist so that 
\beq
\label{36.1}
\beta(\omega_{\tilde\sigma,\tilde y})=0\qquad
\gamma(\omega_{\tilde\sigma,\tilde y})\ge 1/2
\eeq
\beq
\label{36.2}
\beta(\omega_{\bar\sigma,\bar y})=0\qquad
\gamma(\omega_{\bar\sigma,\bar y})= 1/2.
\eeq
\end{lemma}

\no\proof
Set $\forall(\sigma,y)\in\cH$
$$
\theta(\sigma,y)=(\gamma(\omega_{\sigma,y}),\beta(\omega_{\sigma, y}) )
$$
and for all $(\sigma,y)\in\cH$ and $s\in(0,1)$
\beq
\label{36.3}
\cG(\sigma,y,s)=(1-s)(\sigma,y)+s\theta(\sigma,y).
\eeq
Let observe that $\partial\cH$ is isomorphic to a sphere in $\R^4$ and
that $\left({1/ 2},0\right) \in \stackrel{\circ}{\cH}$, so the
line $\left[\left.{1/ 2},+\infty\right)\right.\times\{0\}$ crosses
$\partial\cH$.
Furthermore, (\ref{36.3}) shows $\partial\cH$ is homotopic to
$\theta(\partial\cH)$, thus (\ref{36.1}) is false only if
$(\sigma,y)\in\partial\cH$ and $s\in[0,1]$ exist so that
$\cG(\sigma,y,s)=\left({1/ 2},0\right)$.
On the other hand, in order to prove (\ref{36.2}) it is enough to show
\beq
\label{36.4} 
d(\theta, \stackrel{\circ}{\cH}, (1/2,0))=1,
\eeq
and, since $d(I, \stackrel{\circ}{\cH}, (1/2,0))=1$, (\ref{36.4})
straightly follows by the topological degree homotopy invariance, if
for all $(\sigma,y)\in\partial\cH$, for all $s\in [0,1]$
$\cG(\sigma,y,s)\neq (1/2,0)$.

Therefore, in view of (\ref{14.3}), to prove (\ref{36.1}) and (\ref{36.2}) we just need to
verify
$$
((1-s)\sigma+s\gamma(\omega_{\sigma,y}),(1-s)y+s\beta(\omega_{\sigma,y}))\neq
(1/2,0)\quad\forall (y,\sigma)\in\partial\cH,\ \forall s\in[0,1].
$$
Now $\partial\cH=\cF_1\cup\cF_2\cup\cF_3$ with
$$
\cF_1=\{(\sigma,y)\in\partial\cH\ :\ |y|\le\bar r,\ \sigma=\sigma_1\}
$$
$$
\cF_2=\{(\sigma,y)\in\partial\cH\ :\ |y|\le\bar r,\ \sigma=\sigma_2\}
$$
$$
\cF_3=\{(\sigma,y)\in\partial\cH\ :\ |y|=\bar r,\
\sigma\in[\sigma_1,\sigma_2]\}.
$$
By (\ref{34.1}) we get
$$
\forall (\sigma,y)\in\cF_1\ :\
(1-s)\sigma_1+s\gamma(\omega_{\sigma_1,y})<1/2
$$
$$
\forall (\sigma,y)\in\cF_2\ :\
(1-s)\sigma_2+s\gamma(\omega_{\sigma_2,y})>1/2.
$$

If $(\sigma,y)\in\cF_3$, using (\ref{29.1})$(c)$ we get
$$
((1-s)y+s\beta(\omega_{\sigma,y})|y)=(1-s)|y|^2+s(\beta(\omega_{\sigma,y})|y)>0
$$
hence $(1-s)y+s\beta(\omega_{\sigma,y})\neq 0$.

\qed

\begin{lemma}
\label{L38}
Let $\sigma_1,\sigma_2,\bar r, \cH$ as in Lemma \ref{L34.1}.
Assume that (\ref{stimaprecisabis})  holds, then
\beq
\label{38.1}
L=\sup\{I_0(\hat\omega_{\sigma,y})\ :\
(\sigma,y)\in\cH\}<{2\over 3}S^{3/2}.
\eeq
\end{lemma}

\no\proof
Taking into account that $\omega_{\sigma,y}\in\cN_\infty$ by definition,
$t_{\sigma,y,0}\ge 1$ by (\ref{7.1}),  and using (\ref{23'.1}) we have for
all $(\sigma,y)\in\cH$
\begin{eqnarray}
\nonumber I_0(\hat\omega_{\sigma,y}) & = &
{1\over
  4}\int_{\R^3}[|\D\hat\omega_{\sigma,y}|^2+W(x)\hat\omega^2_{\sigma,y}]dx+{1\over
  12}\int_{\R^3}\hat\omega^6_{\sigma,y}dx
\\
\nonumber
&=& 
{1\over
  4}\int_{\R^3}t^2_{\sigma,y,0}[|\D\omega_{\sigma,y}|^2+W(x)\omega^2_{\sigma,y}]dx+{1\over
  12}t^6_{\sigma,y,0}\int_{\R^3}\omega^6_{\sigma,y}dx
\\ \label{n23.1}
&\le& 
t^6_{\sigma,y,0}\left[{1\over 3}\|\omega\|^2+{1\over
    4}|W|_{3/2}|\omega|^2_6\right]
\\ \nonumber
&\le& 
t^6_{\sigma,y,0}\left[\Sigma+{1\over 4\, S}|W|_{3/2}\|\omega\|^2\right]
\\ \nonumber
&\le& t^6_{\sigma,y,0}\left(1+{3\over 4\, S }|W|_{3/2}\right)\Sigma.
\end{eqnarray}
On the other hand, $\hat\omega_{\sigma,y}\in\cN_0$ implies 
$$
t^2_{\sigma,y,0}\int_{\R^3}[|\D\omega_{\sigma,y}|^2+W(x)\omega^2_{\sigma,y}]dx+
t^4_{\sigma,y,0}\int_{\R^3}K(x)\phi_{\omega_{\sigma,y}}\omega^2_{\sigma,y}dx-
t^6_{\sigma,y,0}\int_{\R^3}\omega^6_{\sigma,y}dx=0
$$
hence, using (\ref{7.1})
$$
\left(\|\omega_{\sigma,y}\|^2+\int_{\R^3}W(x)\omega^2_{\sigma,y}dx+\int_{\R^3}K(x)\phi_{\omega_{\sigma,y}}\omega^2_{\sigma,y}dx\right)-t^2_{\sigma,y,0}|\omega_{\sigma,y}|^6_6\ge
0.
$$
Moreover, since $|\phi_{\omega_{\sigma,y}}|_6\le {1\over S}\,
|K|_2|\omega_{\sigma,y}|^2_6,$ 
\begin{eqnarray}
\nonumber
t^2_{\sigma,y,0}&\le&
1+{|W|_{3/2}\over|\omega|^4_6}+{|K|_2|\phi_{\omega_{\sigma,y}}|_6\over|\omega|_6^4}
\\
\label{n23.2}
&\le& 1+{|W|_{3/2}\over(3\Sigma)^{2/3}}+{1\over S}\, {|K|_2^2\over |\omega|_6^2}
\\ 
\nonumber
&\le &1+{|W|_{3/2}\over S}+{|K|_2^2\over S^{3/2}}.
\end{eqnarray}
So,  by  (\ref{stimaprecisabis}), (\ref{n23.1}),
(\ref{n23.2}) we get
$$
I_0(\hat \omega_{\sigma,y})\le \left(1+{|W|_{3/2}\over S}+{|K|_2^2
 \over   S^{3/2}}\right)^3
\,
\left(1+{3\over 4\, S}|W|_{3/2}\right)\Sigma<{2\over 3}\, S^{3/2}.
$$

\qed

In what follows we use the notation
$$
I_0^c=\{u\in\cN_0\ :\ I_0(u)\le c\}.
$$

{\mbox {{\underline {\sf Proof Theorem \ref{V=0}}}}}
Collecting the results of Proposition \ref{P14.1}, Lemma \ref{L36} and
Lemma \ref{L38}, setting $\hat\omega_{\bar\sigma,\bar y}=t_{\bar\sigma,\bar
  y,0}\omega_{\bar\sigma,\bar y}$, where $\omega_{\bar\sigma,\bar
  y}\in\stackrel{\circ}{\cH}$ is defined in Lemma \ref{L36}, we obtain
$$
{1\over 3}\, S^{3/2}<\bar c< \cB_0\le I_0(\hat\omega_{\bar\sigma,\bar
  y})\le L<{2\over 3}\, S^{3/2}.
$$
Let us show that a critical level of $I_0$ constrained on $\cN_0$
exists in the interval $\left({1\over 3}\, S^{3/2},\right.$ $\left.{2\over 3}\,
    S^{3/2}\right)$. 
We argue by contradiction and we assume there are not critical levels
in $\left({1\over 3}\, S^{3/2},{2\over 3}\, S^{3/2}\right)$.
So, since $I_0$ constrained on $\cN_0$
satisfies the Palais-Smale condition in the energy range $\left({1\over 3}\, S^{3/2},{2\over 3}\,
    S^{3/2}\right)$, by using standard deformation arguments we find a
  number $\delta>0$ such that $\cB_0-\delta>\bar c$, $L+ \delta<{2\over 3}\, S^{3/2}$ and a continuous function
$$
\eta:I_0^{L+\delta}\lo I_0^{\cB_0-\delta}
$$
such that 
$$
\eta(u)=u\qquad\forall u\in I_0^{\cB_0-\delta}.
$$
Then, we remark that 
$$
\forall (\sigma,y)\in\cH\qquad
I_0(\eta(\hat\omega_{\sigma,y}))\le\cB_0-\delta
$$
so
\beq
\label{41.1}
\Theta(\sigma,y):=(\gamma(\eta(\hat\omega_{\sigma,y})),\beta(\eta(\hat\omega_{\sigma,y})))\neq
(0,1/2).
\eeq
On the other hand, by Lemma \ref{L34.1} $\forall
(\sigma,y)\in\partial\cH$
$$
I_0(\hat\omega_{\sigma,y})<\bar c<\cB_0-\delta\quad \Rightarrow\quad
\eta(\hat\omega_{\sigma,y})=\hat\omega_{\sigma,y}
$$
from which
$$
\Theta
(\sigma,y)=\theta(\sigma,y)=(\gamma(\hat\omega_{\sigma,y}),\beta(\hat\omega_{\sigma,y}))\quad\forall
(\sigma,y)\in\partial\cH.
$$
Therefore, by the homotopy invariance of topological degree, we deduce
$$
1=d(\theta,\stackrel{\circ}{\cH},(1/2,0))=d(\Theta,\stackrel{\circ}{\cH},(1/2,0))
$$
that implies the existence of $(\hat\sigma,\hat y)\in\cH$ for which
$$
\Theta(\hat\sigma,\hat y)=(0,1/2)
$$
contradicting (\ref{41.1}).  

To conclude the proof we only must show
that   solutions  corresponding to critical levels lying in the
interval  $\left({1\over 3}\, S^{3/2},{2\over 3}\, S^{3/2}\right)$
cannot change sign. 
Indeed, assume $u =  u^+ - u^-$ is a solution with $u^+ \neq0$, $u^-\neq 0$. 
Then, by Proposition \ref{P7.2} and taking into account that
$\phi_{u^+}\le \phi_u$, we get
\begin{eqnarray}
\label{f1}{1\over 3}\, S^{3/2}
& \le & I(t_{u^+}u^+)\\
\nonumber 
& = & {1\over 3}t_{u^+}^2\, \int_{\R^3}(|\D
u^+|^2+W(x)(u^+)^2)\, dx+{1\over 12}
t_{u^+}^4\int_{\R^3}K(x)\phi_{u^+}(x)(u^+)^2dx \\
\nonumber &\le & 
{1\over 3}t_{u^+}^2\, \int_{\R^3}(|\D
u^+|^2+W(x)(u^+)^2)\, dx+{1\over 12}
t_{u^+}^4\int_{\R^3}K(x)\phi_{u}(x)(u^+)^2dx. 
\end{eqnarray}
We { claim that $t_{u^+}\le 1$}. 
Once the claim is proved, by  (\ref{f1}) we obtain
\beq
\label{f2}
{1\over 3}\, S^{3/2} \le 
{1\over 3}\, \int_{\{u>0\}}(|\D
u|^2+W(x)\, u^2)\, dx+{1\over 12}
\int_{\{u>0\}}K(x)\phi_{u}(x)\, u^2dx.
\eeq
Likewise,
$$
{1\over 3}\, S^{3/2} \le 
{1\over 3}\, \int_{\{u<0\}}(|\D
u|^2+W(x)\, u^2)\, dx+{1\over 12}
\int_{\{u<0\}}K(x)\phi_{u}(x)\, u^2dx,
$$
so that $I(u)\ge {2\over 3}\, S^{3/2}$ follows.

Let us prove that $t_{u^+}\le 1$. 
By definition of projection,
\beq
\label{f4}
t_{u^+}^2\, \int_{\R^3}(|\D
u^+|^2+W(x)(u^+)^2)\,
dx+t_{u^+}^4\int_{\R^3}K(x)\phi_{u^+}(x)(u^+)^2dx 
= t_{u^+}^6 |u^+|_6^6,
\eeq
and, since $u$ is a critical point,
\beq
\label{f5}
|u^+|_6^6= \int_{\R^3}(|\D
u^+|^2+W(x)(u^+)^2)\,
dx+\int_{\R^3}K(x)\phi_{u}(x)(u^+)^2dx .
\eeq
So $t_{u^+}$ is the positive solution of $\vi (t)=0$, where
$$
\vi(t):=t^4|u^+|_6^6-t^2 \int_{\R^3}K(x)\phi_{u^+}(x)(u^+)^2dx +
\int_{\R^3}K(x)\phi_{u}(x)(u^+)^2dx -|u^+|_6^6.
$$
A stright computation shows that $\vi'(t)> 0$ $\forall t\ge 1$,
because, by (\ref{f5}),
\begin{eqnarray*}
{1\over |u^+|_6^6} \int_{\R^3}K(x)\phi_{u^+}(x)(u^+)^2dx &\le&
{1\over |u^+|_6^6} \int_{\R^3}K(x)\phi_{u}(x)(u^+)^2dx \\
& = & 1-{1\over |u^+|_6^6}\int_{\R^3}(|\D
u^+|^2+W(x)(u^+)^2)\,
dx<1.
\end{eqnarray*}
So the claim follows from
$$
\vi(1)=\int_{\R^3}K(x)[\phi_{u}(x)-\phi_{u^+}(x)](u^+)^2dx\ge 0.
$$

\qed

\begin{lemma}
\label{L42}
Let $\bar r,\sigma_1,\sigma_2$ and $\cH$ as in Lemma \ref{L34.1}.
Then there exists a number $\bar V>0$ such that for all $V_\infty\in
(0,\bar V)$
\beq
\label{42.1}
\gamma(\tilde\omega_{\sigma_1,y})
<1/2,\quad\gamma(\tilde\omega_{\sigma_2,y})>1/2,\qquad\forall 
y\in\R^3,\ |y|<\bar r,
\eeq
\beq
\label{42.2}
\tilde l:=\sup\{I(\tilde\omega_{\sigma,y})\ :\ (\sigma,y)\in\partial\cH\}<\bar
c.
\eeq
Furthermore, if (\ref{stimaprecisa})  holds true,
then  $\bar V$ can be found so that, in addition to (\ref{42.1}) and  (\ref{42.2}),
\beq
\label{42.3}
\tilde s:=\sup\{I(\tilde\omega_{\sigma,y})\ :\
(\sigma,y)\in\cH\}<{2\over 3}\, S^{3/2}
\eeq
is satisfied.
\end{lemma}

\no\proof
By (\ref{14.3})
$$
\gamma(\tilde\omega_{\sigma,y})=
\gamma(\hat\omega_{\sigma,y})=
\gamma(\omega_{\sigma,y})\qquad\forall
(\sigma,y)\in(0,+\infty)\times\R^3.
$$
Hence relations (\ref{42.1}) straightly  follow from (\ref{34.1}).
Moreover
$$
1={\|\omega_{\sigma,y}\|^2_{\cD^{1,2}}\over
  |\omega_{\sigma,y}|^6_6}=t^4_{\sigma,y}-{\int_{\R^3}(V_\infty+W(x))\omega^2_{\sigma,y}-t^2_{\sigma,y}\int_{\R^3}K(x)\phi_{\omega_{\sigma,y}}\omega^2_{\sigma,y}\over
  |\omega_{\sigma,y}|^6_6}
$$
and
\beq
\label{32.4}
\int_{\R^3}V_\infty\omega^2_{\sigma,y}(x)\,
dx=V_\infty\sigma^2\int_{B_1(0)}\omega^2(x)\, dx
\eeq
imply 
\beq
\label{43.1}
\lim_{V_\infty\to
  0}\sup_{(\sigma,y)\in\cH}|t_{\sigma,y}-t_{0,\sigma,y}|=0.
\eeq
Then, if $V_\infty$ is suitably small (\ref{42.2}) and (\ref{42.3})
are consequence of (\ref{32.4}), (\ref{43.1}), (\ref{34.2}) and
(\ref{38.1}).

\qed

\no {\mbox {{\underline {\sf Proof of Theorem \ref{T2}}} \hspace{2mm}}}
In what follows $\bar V$ denotes the number whose existence
is stated in Lemma \ref{L42} and we assume $V_\infty \in (0, \bar V).$

To prove the theorem we intend to show that a critical level exists
in the interval $\left({1\over 3}\, S^{3/2},\bar c\right)$ and that if
(\ref{stimaprecisa})  holds another critical level exists in
$\left(\bar c,{2\over 3}\, S^{3/2}\right)$. 

By using (\ref{21.3}), (\ref{36.1}) together with (\ref{14.3}),
(\ref{42.2}), and (\ref{23'.0}) we deduce
\beq
\label{43.2}
{1\over 3}\, S^{3/2}<\mu\le I(\tilde\omega_{\tilde\sigma,\tilde
  y})\le\tilde l<\bar c<\cB_0.
\eeq
Arguing by contradiction, we assume there are no critical levels in  
 $\left({1\over 3}\, S^{3/2},\bar c\right)$.
Since the Palais-Smale compactness condition  holds in that energy
range, we can find a positive number $\delta_1>0$ such that  
$$
\mu-\delta_1>{1\over 3}\, S^{3/2}\qquad\tilde l+\delta_1<\bar c
$$
and a continuous function
$$
\eta:[0,1]\times I^{\tilde l+\delta_1}\lo I^{\tilde l-\delta_1}
$$
such that
\begin{eqnarray}
\nonumber
& \eta(0,u)=u  & \\
\nonumber
& \eta(s,u)=u &\hspace{1cm}\forall u\in I^{\mu-\delta_1}, \ \forall s\in[0,1]
\\
\label{44.1}
&I\circ \eta\, (s,u)\le I(u) &\hspace{1cm} \forall s\in[0,1]
\\
\label{44.2}
& \eta(s,I^{\tilde l+\delta_1})\subset I^{\mu-\delta_1}. &
\end{eqnarray}
Therefore, definition of $\tilde l$ and (\ref{44.2}) give 
\beq
\label{44.3}
(\sigma,y)\in\partial\cH\ \Rightarrow\
I(\tilde\omega_{\sigma,y})\le\tilde l\ \Rightarrow\
I(\eta(1,\tilde\omega_{\sigma,y})\le\mu-\delta_1.
\eeq
Let us consider $\forall s\in [0,1]$, $\forall (\sigma,y)\in\cH$
$$
\Gamma(\sigma,y,s)\ =\ \left\{
\begin{array}{lc}
\cG(\sigma,y,2s) & s\in[0,1/2] \\
\left(\beta\circ\eta\, 
  (2\sigma-1),\tilde\omega_{\sigma,y}),\gamma\circ\eta\, (2\sigma-1,\tilde\omega_{\sigma,y})\right)
& s\in[1/2,1]
\end{array}
\right.
$$
where $\cG$ is the map defined in (\ref{36.3}). 
As already shown in Lemma \ref{L36},
$$
\forall s\in[0,1/2],\ \forall (\sigma,y)\in\partial \cH,\ \Gamma
(\sigma,y,s)\neq (0,1/2).
$$
Furthermore, by using (\ref{43.2}) and (\ref{44.1}) we deduce $\forall
s\in[1/2,1]$ $\forall (\sigma,y)\in\cH$
$$
I(\eta(2s-1,\tilde\omega_{\sigma,y}))\le
I(\tilde\omega_{\sigma,y})\le\tilde l<\bar c<\cB_0\le\cB_{V_\infty}
$$
which gives
$$
\eta(2s-1,\tilde\omega_{\sigma,y})\neq (0,1/2).
$$
Thus, arguing as in Lemma \ref{L36}, we find $(\check \sigma,\check
y)\in\partial\cH$ such that 
$$
\beta\circ\eta\, (1,\tilde \omega_{\check \sigma,\check y})=0,\quad 
\gamma\circ\eta\, (1,\tilde \omega_{\check \sigma,\check y})\ge 1/2.
$$
Then
$$
I( \eta(1,\tilde \omega_{\check \sigma,\check y}))\ge\mu
$$
which contradicts (\ref{44.3}).

Now, let us suppose that (\ref{stimaprecisa})  holds.
By using (\ref{stimaprecisa2}), (\ref{23.1}), (\ref{36.2}) together with
(\ref{14.3}), and (\ref{42.3}) we deduce
\beq
\label{n26.1}
 \bar c<\cB_0\le \cB_{V_\infty}\le I(\tilde \omega_{\bar\sigma,\bar
  y})\le \tilde s<{2\over 3}\, S^{3/2}.
\eeq
 Repeating the argument of the proof of Theorem \ref{V=0},  the
 existence of a critical level lying 
in the energy interval $\left(\bar c,{2\over 3}\,
S^{3/2}\right)$ follows.
The proof is completed showing, exactly as for Theorem \ref{V=0}, that
critical points $u$ such that $I(u) < {2\over 3}\, S^{3/2}$ are
functions that cannot change sign. 

\qed 

{\small {\bf Acknowledgement}. The authors have been supported by the ``Gruppo
Nazionale per l'Analisi Matematica, la Probabilit\`a e le loro
Applicazioni (GNAMPA)'' of the {\em Istituto Nazionale di Alta Matematica
(INdAM)} - Project: Sistemi differenziali ellittici nonlineari
derivanti dallo studio di fenomeni elettromagnetici.

The second author has been supported also by the research project: 
``Consolidate the Foundations: 
Nonlinear Differential Problems and their Applications'', of
the University of Rome ``Tor Vergata''.
}



\end{document}